\documentclass[aos,MSNbibl,seceqn,nameyear,dvips]{arximspdf}
\usepackage{graphicx}

\doi{10.1214/11-AOS942}
\volume{40}
\issue{1}
\pubyear{2012}
\firstpage{45}
\lastpage{72}

\makeatletter

\newtheorem{theorem}{Theorem}[section]

\newtheorem{lem}[theorem]{Lemma}
\newtheorem{prop}[theorem]{Proposition}

\newproclaim{rem}[theorem]{Remark}
\newproclaim{Claim}{Claim}
\newproclaim{Fact}{Fact}

\newcommand{\RR}{\mathbb{R}}

\makeatother

\begin{document}
\begin{frontmatter}

\title{Likelihood based inference for current status data on a grid:
A boundary phenomenon and an adaptive inference procedure}
\runtitle{Asymptotics for current status data}

\begin{aug}
\author[A]{\fnms{Runlong} \snm{Tang}\corref{}\thanksref{t1}\ead[label=e1]{runlongt@princeton.edu}},
\author[B]{\fnms{Moulinath} \snm{Banerjee}\thanksref{t1}\ead[label=e2]{moulib@umich.edu}\ead[label=u2,url]{http://www.stat.lsa.umich.edu/\textasciitilde moulib/}}
\and\\
\author[C]{\fnms{Michael R.} \snm{Kosorok}\ead[label=e3]{kosorok@unc.edu}\ead[label=u3,url]{http://www.bios.unc.edu/\textasciitilde kosorok/}}
\runauthor{R. Tang, M. Banerjee and M. R. Kosorok}
\affiliation{Princeton University, University of Michigan, Ann Arbor, and~University~of~North~Carolina, Chapel Hill}
\address[A]{R. Tang\\
Department of Operations Research \\
\quad and Financial Engineering\\
Princeton University\\
214 Sherrerd Hall, Charlton Street\\
Princeton, New Jersey 08544\\
USA\\
\printead{e1}}
\address[B]{M. Banerjee \\
Department of Statistics \\
University of Michigan, Ann Arbor\\
439 West Hall \\
1085 South University \\
Ann Arbor, Michigan 48109 \\
USA\\
\printead{e2}}
\address[C]{M. R. Kosorok\\
Department of Biostatistics\\
University of North Carolina, Chapel Hill\\
3101 McGavran-Greenberg Hall\\
CB 7420 \\
Chapel Hill, North Carolina 27599 \\
USA\\
\printead{e3}} %adresu isvedimo komanda gale!
\end{aug}

\thankstext{t1}{Supported in part by NSF Grants DMS-07-05288
and DMS-10-07751.}

%and DMS-10-07751.}

% HISTORY:
\received{\smonth{8} \syear{2011}}
\revised{\smonth{11} \syear{2011}}

% ABSTRACT
%
\begin{abstract}
In this paper, we study the nonparametric maximum likelihood estimator
for an event time distribution function at a point in the current
status model with observation times supported on a grid of potentially
unknown sparsity and with multiple subjects sharing the same
observation time. This is of interest since observation time ties occur
frequently with current status data. The grid resolution is specified
as $c n^{-\gamma}$ with $c > 0$ being a scaling constant and $\gamma>0$
regulating the sparsity of the grid relative to $n$, the number of
subjects. The asymptotic behavior falls into three cases depending on
$\gamma$: regular Gaussian-type asymptotics obtain for $\gamma< 1/3$,
nonstandard cube-root asymptotics prevail when $\gamma> 1/3$ and
$\gamma=1/3$ serves as a~boundary at which the transition happens. The
limit distribution at the boundary is different from either of the
previous cases and converges weakly to those obtained with
$\gamma\in(0,1/3)$ and $\gamma\in(1/3,\infty)$ as~$c$ goes to
$\infty$
and $0$, respectively. This weak convergence allows us to develop an
adaptive procedure to construct confidence intervals for the value of
the event time distribution at a point of interest \textit{without
needing to know or estimate $\gamma$}, which is of enormous advantage
from the perspective of inference. A simulation study of the adaptive
procedure is presented.
\end{abstract}

% KEYWORDS
%

\end{frontmatter}

%s1 #&#
\section{Introduction}\label{intro}
The current status model is one of the most well-studied survival
models in statistics. An individual at risk for an event of interest is
monitored at a random observation time, and an indicator of whether the
event has occurred is recorded. An interesting feature of this kind of
data is that the underlying event time distribution, $F$, can be
estimated by its nonparametric maximum likelihood estimator (NPMLE) at
only $n^{1/3}$ rate when the observation time is a continuous random
variable. Under mild conditions on $F$, the limiting distribution of
the NPMLE in this setting is the non-Gaussian Chernoff distribution:
the distribution of the unique minimizer of $\{W(t) + t^2\dvtx t
\in\mathbb{R}\}$, where $W(t)$ is standard two-sided Brownian motion.
This is in contrast to data with right-censored event times where $F$
can be estimated nonparametrically at rate $\sqrt{n}$ and is ``pathwise
norm-differentiable'' in the sense of \citet{Vaart1991}, admitting
regular estimators and normal limits. Interestingly, when the
observation time distribution has finite support, the NPMLE for $F$ at
a point asymptotically simplifies to a binomial random variable and is
also $\sqrt{n}$ estimable and regular, with a normal limiting
distribution.

An extensive amount of work has been done for inference in the current
status model under the assumption of a continuous distribution for the
observation time: the classical model considers $n$ subjects whose
survival times $T_1, T_2, \ldots, T_n$ are i.i.d. $F$ and whose
inspection times
$X_1, X_2, \ldots, X_n$ are i.i.d. with some continuous distribution,
say $G$; furthermore, in the absence of covariates, the $X_i$'s and
$T_i$'s are
considered mutually independent. The observed data are $\{\Delta_i,
X_i\}_{i=1}^n$, where $\Delta_i = 1\{T_i \leq X_i\}$, and one is
interested in estimating $F$ as $n$ goes to infinity. More
specifically, for inference on the value of $F$ at a pre-fixed point of
interest under a continuous observation time, see, for example,
\citet{Groeneboom1992}, who establish the convergence of the normalized NPMLE
to Chernoff's distribution; \citet{Keiding1996}; \citet{Wellner2000},
who develop pseudo-likelihood estimates of the mean function of a
counting process with panel count data, current status data being a
special case; \citet{Banerjee2001} and \citet{Banerjee2005}, who
develop an asymptotically pivotal likelihood ratio based method; \citet
{Sen2006}, who extend the results of \citet{Wellner2000} to
asymptotically pivotal inference for $F$ with mixed-case
interval-censoring; and \citet{Groeneboom2010} for smoothed isotonic
estimation, to name a few.

However, somewhat surprisingly, the problem of making inference on $F$
when the observation times lie on a grid with multiple subjects sharing
the same observation time has never been satisfactorily addressed in
this rather large literature. This important scenario, which transpires
when the inspection times for individuals at risk are evenly spaced,
and multiple subjects can be inspected at any inspection time, is
\textit{completely precluded} by the assumption of a continuous $G$, as this
does not allow ties among observation times. Consider, for example, a
tumorigenicity study where a large number of mice are exposed to some
carcinogen at a particular time, and interest centers on the time to
development of a tumor. A typical procedure here would be to randomize
the mice to be sacrificed over a number of days following exposure; so,
one can envisage a protocol of sacrificing a fixed number~$m$ of mice
at 24 hrs post-exposure, another $m$ mice at 48 hours and so on. The
sacrificed mice are then dissected and examined for tumors, thereby
leading to current status data on a grid. A pertinent question in this
setting is: what is the probability that a mouse develops a tumor by an
$M$-day period after exposure? This involves estimating $F(24 M)$,
where $F$ is the distribution function of the time to
tumor-development. Similar grid-based data can occur with human
subjects in clinical settings.

In this paper we provide a clean solution to this problem based on the
NPMLE of $F$ which, as is well known, is obtained through isotonic regression
[see, e.g., \citet{Robertson1988}]. The NPMLE of $F$ in the current
status model (and more generally in nonparametric monotone function
models) has a long history and has been studied extensively. In
addition to the attractive feature that it can be computed without
specifying a bandwidth, the NPMLE of $F(x_0)$ (where $x_0$ is a fixed
point) attains the best possible convergence rate, namely $n^{1/3}$, in
the ``classical'' current status model with continuous observation
times, under the rather mild assumption that $F$ is continuously
differentiable in a neighborhood of $x_0$ and has a nonvanishing
derivative at $x_0$. \textit{This rate cannot be bettered by a smooth
estimate under the assumption of a single derivative.} As demonstrated
in \citet{Groeneboom2010}, smoothed monotone estimates of $F$ can
achieve a faster $n^{2/5}$ rate under a \textit{twice-differentiability
assumption on $F$}; hence, the faster rate requires additional
smoothness. However, as we wish to approach our problem under minimal
smoothness assumptions, the isotonic NPMLE is the more natural choice.
(Smoothing the NPMLE would introduce an exogenous tuning parameter
without providing any benefit from the point of view of the convergence rate.)

The key step, then, is to determine the best asymptotic approximation
to use for the NPMLE in the grid-based setting discussed above. If, for
example, the number of observation times, $K$, is far smaller than $n$,
the number of subjects, the problem is essentially a parametric one, and
it is reasonable to expect that normal approximations to the MLE will
work well. On the other hand, if $K = n$, that is, we have a very fine
grid with each
subject having their own inspection time, the scenario is similar to
the current status model with continuous observation times where no two
inspection times coincide, and one may expect a Chernoff approximation
to be adequate. However, there is an entire spectrum of situations in
between these extremes depending on the size of the grid,~$K$, relative
to~$n$, and if~$n$ is ``neither too large, nor too small relative to~$K$,''
neither of these two approximations would be reliable.

Some work on the current status model or closely related variants under
discrete observation time settings should be noted in this context.
\citet{Yu1998} have studied the asymptotic properties of the NPMLE
of $F$ in the current status model with discrete observation times, and
more recently \citet{Maathuis2010} have considered nonparametric
inference for (finitely many) competing risks current status data under
discrete or grouped observation times. However, these papers consider
situations where the observation times are i.i.d. copies from a
\textit{fixed discrete distribution} (but not necessarily finitely
supported) on the time-domain and are therefore not geared toward
studying the effect of the trade-off between $n$ and~$K$, that is,
\textit{the effect of the relative sparsity of the number of distinct
observation times to the size of the cohort of individuals} on
inference for $F$. In both these papers, the pointwise estimates of $F$
are $\sqrt{n}$ consistent and asymptotically normal; but as \citet
{Maathuis2010} demonstrate in Section 5.1 of their paper, when the
number of distinct observation times is large relative to the sample
size, the normal approximations are suspect.

Our approach is to couch the problem in an asymptotic framework
where~$K$ is allowed to increase with $n$ at rate $n^{\gamma}$ for
some $0 < \gamma\leq1$ and study the behavior of the NPMLE at a
grid-point. This is achieved by considering the current status model on a
regular grid over a compact time interval, say $[a,b]$, with unit
spacing $\delta\equiv\delta_n = c n^{-\gamma}$, $c$ being a scale
parameter.
It will be seen that the limit behavior of the NPMLE depends heavily on
the ``sparsity parameter'' $\gamma$, with the Gaussian approximation
prevailing for $\gamma< 1/3$ and the Chernoff approximation for
$\gamma> 1/3$. When $\gamma= 1/3$, one obtains a~discrete analog of
the Chernoff distribution which depends on $c$. Thus, there is an
entire family of what we call \textit{boundary distributions}, indexed
by $c$, say $\{F_c\dvtx c > 0\}$, by manipulating which, one can approach
either the Gaussian or the Chernoff. As $c$ approaches 0, $F_c$
approximates the Chernoff while, as
$c$ approaches $\infty$, it approaches the Gaussian. This property
allows us to develop an \textit{adaptive procedure} for setting confidence
intervals for the value of $F$ at a grid-point that \textit{obviates the
need to know or estimate $\gamma$}, the critical parameter in this
entire business as it completely dictates the ensuing asymptotics. The
adaptive procedure involves \textit{pretending} that the true unknown
underlying unknown $\gamma$ is at the boundary value $1/3$, computing
a~\textit{surrogate} $c$, say~$\hat{c}$, by equating $(b-a)/K$, the
spacing of the grid (which is
computable from the data), to $\hat{c} n^{-1/3}$ and using~$F_{\hat
{c}}$, to approximate the distribution of the appropriately normalized
NPLME. The details are given in Section \ref{sec4}. It is seen that this
procedure provides asymptotically correct confidence intervals
regardless of the true value of $\gamma$. Our procedure does involve
estimating some nuisance parameters, but this is readily achieved via
standard methods.

The rest of the paper is organized as follows. In Section \ref{sec2}, we present
the mathematical formulation of the
problem and introduce some key notions and characterizations. Section
\ref{sec3}
presents the main asymptotic results and
their connections to existing work. Section \ref{sec4} addresses the important
question of
\textit{adaptive inference} in the current status model: given a
time-domain and current status data observed at times on a
regular grid of an \textit{unknown level of sparsity} over the domain,
how do we
make inference on $F$? Section \ref{sec5} discusses the implementation of the
procedure and presents results from simulation studies, and Section
\ref{sec6}
concludes with a discussion of the findings of this paper and their
implications for monotone regression models in general, as well as more
complex forms
of interval censoring and interval censoring with competing risks. The
\hyperref[app]{Appendix} contains some technical
details.
% \notes
% \begin{itemize}
% \item Perhaps we also need to include some references about the
% grouped data problems.
% \end{itemize}

% Chernoff's distribution is the same as that of half the left
%derivative of the convex minorant of $X(t)$ at 0 which we shall
%henceforth
%denote as $\mbox{slo-gcm}(X)(0)$, slo-gcm standing for `slope of
%convex minorant'.

%s2 #&#
\section{Formulation of the problem}\label{sec2}

Let $\{T_{i,n}\}_{i=1}^{n}$ be i.i.d. survival times
following some unknown
distribution $F$ with Lebesgue density $f$ concentrated
on the time-domain
$[a',b']$ with
$0 \leq a' < b' < \infty$ (or supported on $[a', \infty)$ if no such
$b'$ exists) and
$\{X_{i,n}\}$ be i.i.d. observation times
drawn from a discrete probability measure $H_{n}$ supported on a regular
grid on $[a,b]$ with $a'\leq a < b < b'$.
Also, $T_{i,n}$ and $X_{i,n}$ are assumed to be
independent for each $i$.
However, $\{T_{i,n}\}$ are not
observed; rather, we observe $\{Y_{i,n} = 1\{T_{i,n} \leq X_{i,n} \} \}$.
This puts us in the setting of a binary regression model
with $Y_{i,n} | X_{i,n} \sim\operatorname{Bernoulli}(F(X_{i,n}))$.
% with $F$ increasing.
We denote the support of $H_n$ by $\{t_{i,n}\}_{i=1}^K$ where
the $i$th grid point
$t_{i,n} = a+i\delta$,
the unit spacing
$\delta= \delta(n) = c n^{-\gamma}$
(also referred to as the \textit{grid resolution})
with $\gamma\in(0,1]$ and $c > 0$,
and the number of grid points $K = K(n) =\lfloor(b-a)/\delta\rfloor$.
On this grid, the distribution $H_{n}$
is viewed as a discretization of an absolutely continuous distribution $G$,
whose support contains $[a,b]$
and whose Lebesgue density is denoted as $g$.
More specifically, $H_{n}\{t_{i,n}\} =
G(t_{i,n}) - G(t_{i-1,n})$, for $i = 2,3,\ldots,K-1$,
$H_{n}\{t_{1,n}\}=G(t_{1,n})$ and $H_{n}\{t_{K,n}\} = 1 - G(t_{K-1,n})$.
For simplicity, these discrete probabilities
are denoted as $p_{i,n}=H_{n}\{t_{i,n}\}$ for each~$i$.
%Suppose $x_{0}\in(a,b)$ is a point
%around which we are interested in determining the
%properties of the MLE of $F$ as a function of the \emph{resolution
%of the grid}, which under our formulation is of the order
%$n^{-\gamma}$, $\gamma$ being the resolution parameter.
%Let $x_{0}\in(a,b)$ be a point
%around which we are interested in determining the
%properties of the NPMLE of $F$
%as a function of the grid resolution $\delta$ or equivalently in terms
%of
%$(c,\gamma)$.
In what
follows, we refer to the pair $(X_{i,n}, Y_{i,n})$ as $(X_i, Y_i)$,
suppressing the dependence on $n$,
%for the sake of notational simplicity
but the \textit{triangular array nature} of our observed
data should be kept in mind.
Similarly, the subscript $n$ is suppressed elsewhere when no confusion
will be caused.
%We also denote, for clarification, the \emph{true unknown survival
%distribution} by $F_{0}$
%with the corresponding Lebesgue density $f_{0}$,
%while the symbol $F$ and $f$ will be used generically.

Our interest lies in estimating $F$ at a grid-point. Since we allow the
grid to change with $n$, this will be accomplished by
specifying a grid-point with respect to a fixed time $x_0 \in(a,b)$
which does not depend on $n$ and can be viewed as an
``anchor-point.'' Define $t_{l} = t_{l,n}$ to be the largest grid-point
less than or equal to~$x_{0}$. We devote our interest to $\hat
{F}(t_l)$. More specifically, we are interested in the limit
distribution of $\hat{F}(t_{l}) - F(t_{l})$ under appropriate normalization?
To this end, we start with the characterization of the NPMLE in this
model. While this is well known from the current status literature, we
include a description tailored for the setting of this paper.

The likelihood function of the data $\{(X_{i}, Y_{i})\}$ is
given by
\[
L_{n}(F)= \prod_{j=1}^{n} F(X_{j})^{Y_{j}}\bigl(1-F(X_{j})\bigr)^{1-Y_{j}}p_{\{
i\dvtx X_{j}=t_{i}\}}
=\prod_{i=1}^{K} F_{i}^{Z_{i}}(1-F_{i})^{N_{i}-Z_{i}}p_{i}^{N_{i}},
\]
where $p_{\{i\dvtx X_{j}=t_{i}\}}$ denotes
the probability that $X_{j}$
equals a genetic grid point~$t_{i}$,
$F_{i}$ is an abbreviation for $F(t_{i})$,
$N_{i}=\sum_{j=1}^{n}\{X_{j}=t_{i}\}$ is the number of observations at $t_{i}$,
$Z_{i}=\sum_{j=1}^{n}Y_{j}\{X_{j}=t_{i}\}$ is
the sum of the responses at $t_{i}$,
$\{\cdot\}$ stands for both a set and its indicator function
with
the meaning depending
on the context and $F$ is generically understood as either a distribution
or the vector $(F_{1}, F_{2}, \ldots, F_{K})$,
which sometimes is also written as $\{F_{i}\}_{i=1}^{K}$.
Then, the log-likelihood function is given by
\[
l_{n}(F) = \log(L_{n}(F)) = \sum_{i=1}^{K}N_{i}\log p_{i} +
\sum_{i=1}^{K}\{[\bar Z_{i}\log F_{i} + (1-\bar Z_{i}) \log
(1-F_{i})^{}]N_{i}\},
\]
where $\bar Z_{i}=Z_{i}/N_{i}$ is the average of the responses at $t_{i}$.

Denote the basic shape-restricted maximizer as
\[
\{F^{\star}_{i}\}_{i=1}^{K} =
\mathop{\arg\max}_{F_{1}\leq\cdots\leq F_{K}}
l_{n}(F).
\]
From the theory of isotonic regression [see, e.g., \citet{Robertson1988}],
we have
\[
\mathop{\arg\max}_{F_{1}\leq\cdots\leq F_{K}}l_{n}(F)
% = \underset{f_{1}\leq\cdots\leq f_{K}}{\mbox{argmax}}l_{n}(f)
= \mathop{\arg\min}_{F_{1}\leq\cdots\leq F_{K}}\sum
_{i=1}^{K}[(\bar Z_{i} - F_{i})^{2}N_{i}].
\]
Thus, $\{F^{\star}_{i}\}_{i=1}^{K}$ is the weighted isotonic
regression of $\{\bar{Z}_i\}_{i=1}^K$ with weights $\{N_i\}_{i=1}^K$,
and exists uniquely. We conventionally define the shape-restricted
NPMLE of $F$ on $[a,b]$ as the following
right-continuous step function:
%
%e2.1 #&#
\begin{equation}\label{equhat-F}
\hat F(t) = \cases{
0, &\quad if $t\in[a,t_{1})$;\cr
F^{\star}_{i}, &\quad if $t\in[t_{i},t_{i+1})$, $i=1,\ldots,K-1$;\cr
F^{\star}_{K}, &\quad if $t \in[t_{K},b]$.}
\end{equation}
Next, we provide a characterization of $\hat F$ as the slope of
\textit{the greatest convex minorant} (GCM)
of a random processes, which proves
useful for deriving the asymptotics for $\gamma\in[1/3,1]$.
Define, for $t\in[a,b]$,
%
%e2.2 #&#
\begin{equation} \label{equGn-Vn}
G_{n}(t)=\mathbb{P}_{n}\{x\leq t\},\qquad V_{n}(t)=\mathbb{P}_{n}y\{
x\leq t\},
\end{equation}
where $\mathbb{P}_{n}$ is the empirical probability measure
based on the data $\{(X_{i}, Y_{i})\}$.
Then, we have, for each $x\in[a,b]$,
%
%e2.3 #&#
\begin{equation} \label{equhat-FLSGCM}
\hat F(x)
= \operatorname{LS}\bigl[\operatorname{GCM} \{(G_{n}(t),V_{n}(t)),t\in[a,b] \}\bigr]
(G_{n}(x)) .
\end{equation}
In the above display, GCM means
the greatest convex minorant
of a set of
points in $\mathbb{R}^{2}$.
For any finite collection of points in $\RR^2$,
its
GCM is a~continuous piecewise linear convex function, and $\operatorname{LS}[\cdot]$
denotes the \textit{left slope or derivative function} of
a convex function.
The term GCM will also be used to refer to the greatest convex minorant
of a real-valued function defined on a sub-interval of the real
line.\vadjust{\goodbreak}
%For such a function $H$, $GCM(H)$ will denote the greatest convex
%minorant of $H$.

Finally, we introduce a number of random processes that will appear in
the asymptotic descriptions of $\hat{F}$.

For constants $\kappa_{1}>0$ and $\kappa_{2}>0$, denote
%
%e2.4 #&#
\begin{equation}
\label{cont-time-proc}
X_{\kappa_{1},\kappa_{2}}(h) =
\kappa_{1} W(h)+ \kappa_{2} h^{2}\qquad \mbox{for } h\in\mathbb{R}  ,
\end{equation}
where $W$ is a two-sided Brownian motion with $W(0)=0$.
Let
$G_{\kappa_{1}, \kappa_{2}}$ be the GCM of
$X_{\kappa_{1},\kappa_{2}}$. Define, for $h\in\mathbb{R}$,
%
%e2.5 #&#
\begin{equation}\label{cont-time-slope}
g_{\kappa_{1},\kappa_{2}}(h) = \operatorname{LS}[G_{\kappa_{1}, \kappa_{2}}](h) .
\end{equation}
%
%$g_{\kappa_{1},\kappa_{2}}(h)=LS(G_{\kappa_{1}, \kappa_{2}})(h)$
%and
%$g_{\kappa_{1},\kappa_{2}}^{o}(h)=
%(\operatorname{LS}(G_{\kappa_{1},\kappa_{2}}^{l})(h)\wedge0)\{h\in(-\infty,0)\} +
%(\operatorname{LS}(G_{\kappa_{1},\kappa_{2}}^{r})(h)\vee0)\{h\in(0,\infty)\}$
%for $h\in\mathbb{R}$.
% Perhaps we can provide several characterizations
% for $\{F^{\star}_{i}\}$.
The process $g_{\kappa_1, \kappa_2}$ will characterize the asymptotic
behavior of a localized NPMLE process in the vicinity of $t_l$ for
$\gamma> 1/3$, from which the large sample distribution of $\hat
{F}(t_l)$ can be deduced.

We also define a three parameter family of processes in discrete time
which serve as discrete versions of the continuous-time processes above.
For $c, \kappa_1$, $\kappa_2 > 0$, let
%
%e2.6 #&#
\begin{eqnarray}
\label{discrete-time-process}
\mathcal{P}_{c,\kappa_1,\kappa_2}(k) &=& (\mathcal{P}_{1,c,\kappa
_1,\kappa_2}(k), \mathcal{P}_{2,c,\kappa_1,\kappa_2}(k)) \nonumber\\[-8pt]\\[-8pt]
&=&
\{ck, \kappa_1 W(ck) + \kappa_2 c^2 k(1+k)\}_{k \in\mathbb{Z}}.
\nonumber
\end{eqnarray}
Define
%
%e2.7 #&#
\begin{equation}
\label{slope-discrete}
\mathbb{X}_{c,\kappa_1,\kappa_2}(ci) = \operatorname{LS}[\operatorname{GCM} \{\mathcal
{P}_{c,\kappa_1,\kappa_2}(k)\dvtx k \in\mathbb{Z}\}](ci)  .
\end{equation}
This slope process will characterize the asymptotic behavior of the
NPMLE in the case $\gamma= 1/3$.

%%\newpage
%
%%\newpage

%s3 #&#
\section{Asymptotic results}\label{sec3}

In this section, we state and discuss results on the asymptotic
behavior of $\hat{F}(t_l)$ for $\gamma$ varying in $(0,1]$. In all
that follows, we make the \textit{blanket assumption}
that $F$ is once continuously differentiable in a neighborhood of $x_{0}$.

%s3.1 #&#
\subsection{\texorpdfstring{The case $\gamma< 1/3$}{The case gamma < 1/3}}\label{secasymptoticsgammaless033}
We start with some technical assumptions:
%requires that $f$ and $g$ have a lower bound `globally'.
%The assumption (\textbf{A1.3}) is purely technical
%and can be replaced by another assumption
%which allow $a'=a$ but require that $F(t_{1})$ does not go
%to 0 too fast. }
%
\begin{longlist}[(A1.3)]
\item[(A1.1)] $F$ has a bounded density $f$ on $[a,b]$, and there
exists $f_{l}>0$
such that $f(x)>f_{l}$ for every $x\in[a,b]$.
\item[(A1.2)] $G$ has a bounded density $g$ on $[a,b]$, and there
exists $g_{l}>0$
such that $g(x)\geq g_{l}$ for every $x\in[a,b]$.
\item[(A1.3)] $a'<a$ and $F(a)>0$.
\end{longlist}
The above assumptions are referred to collectively as (A1).
Letting $t_r$ denote the first grid-point to the right of $t_l$, we
have the following theorem.
%
%th3.1 #&#
\begin{theorem}\label{thmcase1hat-Flimiting-distribution}
If $\gamma\in(0,1/3)$ and \textup{(A1)} holds,
\[
\bigl(\sqrt{N_{l}}\bigl(\hat F(t_{l}) - F(t_{l})\bigr), \sqrt{N_{r}}\bigl(\hat
F(t_{r}) - F(t_{r})\bigr)\bigr)\stackrel{d}{\rightarrow} \sqrt
{F(x_{0})\bigl(1-F(x_{0})\bigr)}
N(0, I_{2}),
\]
where $I_{2}$ is the $2 \times2$ identity matrix.
\end{theorem}

The proof of this theorem is provided in the supplement to this paper
[\citet{Tang2011}]. However, a number of remarks in connection with the
above theorem are in order.
%
%re3.2 #&#
\begin{rem}\label{rem3.2}
From Theorem \ref{thmcase1hat-Flimiting-distribution},
the quantities $\hat{F}(t_{l})$ and $\hat{F}(t_{r})$
with proper centering and scaling
are asymptotically uncorrelated
and independent. In fact,
they are essentially the averages of the responses at the two grid points
$t_{l}$ and $t_{r}$
and are therefore based on responses corresponding to different sets of
individuals.
Consequently, there is no dependence between them in the long run.
Intuitively speaking, $\gamma\in(0, 1/3)$
corresponds to very sparse grids
with successive grid points \textit{far enough}
so that the responses at different grid points fail to influence each other.

It can be shown that for $\gamma\in(0, 1/3)$, $N_{l}/(np_{l})$
converges to 1 in probability
and that $np_{l}/cg(x_{0})n^{1-\gamma}$ converges to 1.
Then the result of Theorem \ref{thmcase1hat-Flimiting-distribution}
can be rewritten as follows:
%(
%n^{\frac{1-\gamma}{2}}(\hat F(t_{l}) - F(t_{l})),
%n^{\frac{1-\gamma}{2}} (\hat
%F(t_{r}) - F(t_{r})))
%
%e3.1 #&#
\begin{equation}\label{equcase1hat-Flimiting-distribution2nd-formation}
\bigl(
n^{(1-\gamma)/2}\bigl(\hat F(t_{l}) - F(t_{l})\bigr),
n^{(1-\gamma)/2} \bigl(\hat
F(t_{r}) - F(t_{r})\bigr)\bigr)
\stackrel{d}{\rightarrow} \alpha c^{-{1}/{2}}N(0, I_{2}),\hspace*{-25pt}
\end{equation}
where $\alpha= \sqrt{F(x_{0}) (1-F(x_{0}))/g(x_{0})}$.
This formulation will be used later, and the parameter $\alpha$ will
be seen to play a critical role in the asymptotic behavior of $\hat
{F}(t_l)$ when $\gamma\in[1/3,1]$ as well.
\end{rem}
%
%re3.3 #&#
\begin{rem}
The proof of the above theorem relies heavily on the below proposition
which deals with the vector of
average responses at the the grid-points: $\{\bar{Z}_i\}
_{i=1}^k$. Since $\bar{Z}_i$ is not
defined when $N_i = 0$, to avoid ambiguity we set $\bar{Z}_i = 0$
whenever this happens. This can be done without affecting the asymptotic
results, since it can be shown that the probability of the event
$\{N_i > 0, i = 1,2,\ldots,K\}$ goes to 1.
%
%pr3.4 #&#
\begin{prop}\label{BRGDYis}
If $\gamma\in(0,1/3)$ and \textup{(A1)} holds, we have
\[
P(\bar Z_{1} \leq\bar Z_{2} \leq\cdots\leq\bar Z_{K} )
\rightarrow1.
\]
\end{prop}

This proposition is established in the supplement, \citet{Tang2011}. It
says that with probability going to 1, the vector $\{\bar{Z}_i\}
_{i=1}^k$ is ordered, and therefore the isotonization algorithm
involved in finding the
NPMLE
of $F$ yields
$\{F_i^{\star}\}_{i=1}^K = \{\bar{Z}_i\}_{i=1}^K$
with probability going to 1. In other words,
asymptotically, isotonization has no effect, and the naive estimates
obtained by averaging the responses at each grid point produce the
NPMLE.
This lemma is really at the heart of the asymptotic derivations for
$\gamma< 1/3$ because it effectively reduces the problem of studying
the $F_i^{\star}$'s, which are obtained through a complex nonlinear
algorithm, to the study of the asymptotics of the $\bar{Z}_i$,
which are linear statistics and can be handled readily using standard
central limit theory. A phenomenon, similar to the one in the above\vadjust{\goodbreak}
proposition, was observed by \citet{Kiefer1976} in connection with
estimating the magnitude of the difference between the empirical
distribution function and its least concave
majorant for an i.i.d. sample from a~concave distribution function. See
Theorem 1 of their paper and the preceding Lemma~4, which establish the
concavity of a piecewise linear estimate of the true distribution
obtained by linearly interpolating the restriction of the empirical
distribution to a grid with spacings of order slightly larger
than~$n^{-1/3}$, $n$ being the sample size. A similar result was obtained in
Lemma 3.1 of \citet{Zhang2001} in connection with isotonic estimation
of a decreasing density when the exact observations are not available;
rather, the numbers of data-points that fall into equi-spaced bins are observed.
\end{rem}

%s3.2 #&#
\subsection{\texorpdfstring{The case $\gamma\in(1/3, 1]$}{The case gamma in (1/3, 1]}}
\label{secasymptoticsgammalarger033}
Our treatment will be condensed since the asymptotics for this case
follow the same patterns as
when the observation times possess a Lebesgue density. That this ought
to be the case is suggested, for example,
by Theorem 1 in \citet{Wright1981}; see, in particular, the condition
on the rate of convergence of the empirical distribution
function of the regressors to the true distribution function in the
case that $\alpha= 1$ in that theorem, which corresponds to
the setting $\gamma> 1/3$ in our problem. Note that the $\alpha$ in
the previous sentence refers to notation in \citet{Wright1981} and
\textit{should not be confused with} the $\alpha$ defined in this
paper.\vspace*{2pt}

In order to study the asymptotics of
the isotonic regression estimator~$\hat{F}(t_l)$, the following
localized process will be of interest:
for
$u \in I_{n}=[(a-t_{l})n^{1/3},\allowbreak(b-t_{l})n^{1/3}]$,
define
%
%e3.2 #&#
\begin{equation}\label{equgamma-0-033Xn}
\mathbb{X}_{n}(u) = n^{1/3}\bigl(\hat F(t_{l}+un^{-1/3})-F(t_{l})\bigr).
\end{equation}
Next, define the following normalized processes on $I_{n}$:
%
%e3.3 #&#
%e3.4 #&#
\begin{eqnarray}
\label{equGn-star}
G_{n}^{\star}(h)
&=& g(x_{0})^{-1}n^{1/3} \bigl(G_{n}(t_l + h n^{-1/3}) - G_{n}(t_l)\bigr),
\\
\label{equVn-star}
V_{n}^{\star}(h)
&=&g(x_{0})^{-1}n^{2/3}
\bigl[
V_{n}(t_l + h n^{-1/3}) - V_{n}(t_l) \nonumber\\[-8pt]\\[-8pt]
&&\hspace*{56.7pt}{} - F(t_{l}) \bigl(G_{n}(t_l + h n^{-1/3}) - G_{n}(t_l)\bigr)
\bigr] \nonumber.
\end{eqnarray}
After some straightforward algebra,
from (\ref{equhat-FLSGCM}) and
(\ref{equgamma-0-033Xn}), we have the following technically useful
characterization
of $\mathbb{X}_{n}$:
for $u\in I_{n}$,
%( G_{n}^{\star}(h), V_{n}^{\star}(h) ),
%h\in I_{n}\} (G_{n}^{\star}(u)) .
%
%e3.5 #&#
\begin{equation}\label{formulaXn}
\mathbb{X}_{n}(u) =
\operatorname{LS} [\operatorname{GCM}
( G_{n}^{\star}(h), V_{n}^{\star}(h) ),
h\in I_{n} ] (G_{n}^{\star}(u)).
\end{equation}
Let $\alpha$ be defined as Remark \ref{rem3.2} and $\beta= f(x_0)/2$.
%and let
%$\beta=f(x_{0})/2$. Let `$\rightsquigarrow$' denote `weak convergence'
%in addition to $\stackrel{d}{\rightarrow}$.
We have the following theorem on the distributional convergence of
$\mathbb{X}_n$.
%
%th3.5 #&#
\begin{theorem}[(Weak convergence of $\mathbb{X}_{n}$)]
\label{thmgamma-033-1Xn-Yn}
Suppose $F$ and $G$ are continuously
differentiable in a neighborhood of $x_{0}$ with derivatives $f$ and $g$.
Assume that $f(x_0) > 0, g(x_0) > 0$ and that $g$ is Lipschitz
continuous in a~neighborhood of $x_0$. Then, the finite-dimensional
marginals of the process~$\mathbb{X}_{n}$
converge weakly to those of the process $g_{\alpha,\beta}$.\vadjust{\goodbreak}
%where $\alpha=\sqrt{F(x_{0})(1-F(x_{0}))/g(x_{0})}$,
%$\beta=f(x_{0})/2$.
\end{theorem}
%
%re3.6 #&#
\begin{rem}
\label{remgamma-greaterThan-oneThirdspecial}
Note that $\mathbb{X}_n(0) = n^{1/3}(\hat{F}(t_{l}) - F(t_{l}))$.
By Theorem \ref{thmgamma-033-1Xn-Yn},
it converges in distribution to $g_{\alpha,\beta}(0)$.
By the Brownian scaling results
on page 1724 of \citet{Banerjee2001}, for $h\in\mathbb{R}$,
\[
g_{\alpha,\beta}(h)
\stackrel{d}{=}
(\alpha^2 \beta)^{1/3} g_{1,1}\bigl((\beta/\alpha)^{2/3}h\bigr)  .
\]
Then, by noting that $g_{1,1}(0)\stackrel{d}{=}2\mathcal{Z}$,
we have the following result:
%
%e3.6 #&#
\begin{equation}
\label{resultgamma-033-1Xn0}
n^{1/3}\bigl(\hat{F}(t_l) - F(t_l)\bigr)
\stackrel{d}{=}
\biggl(\frac{4f(x_{0})F(x_{0})(1-F(x_{0}))}{g(x_{0})}\biggr)^{1/3}
\mathcal{Z}.
\end{equation}
Thus, the limit distribution of $\hat F(t_{l})$
%,the NPMLE of $F(t_{l})$,
is exactly the same as one would encounter in
the
current status model with survival distribution $F$ and
the observation times drawn from a Lebesgue density
function $g$. The proof of this theorem is omitted
as it can be established via arguments similar to those in \citet{Banerjee2007d}
using continuous mapping theorems for slopes of greatest convex minorants.
\end{rem}

%s3.3 #&#
\subsection{\texorpdfstring{The case $\gamma= 1/3$}{The case gamma = 1/3}}\label{secasymptoticsgammaequal033}

Now, we consider the most interesting boundary case $\gamma= 1/3$.
Let the localized process $\mathbb{X}_n(u)$ be defined exactly as in
the previous subsection.
The order of the grid-spacing $\delta$ is now exactly $n^{-1/3}$,
which is the order of localization
around $t_l$ used to define the process $\mathbb{X}_n$, and it follows
that $\mathbb{X}_{n}$ has potential jumps only at
$ci$ for $i\in\mathcal{I}_{n}=(I_{n}/c)\cap\mathbb{Z}$, and it
suffices to consider $\mathbb{X}_{n}$ on those $ci$'s.
For $i \in\mathcal{I}_{n}$,
%
%e3.7 #&#
%e3.8 #&#
\begin{eqnarray}\label{equcase-threeXn}
\mathbb{X}_{n}(ci) &=& n^{1/3}\bigl(\hat{F}(t_l + ci n^{-1/3}) - F(t_l)\bigr)
\\
&=&
\operatorname{LS} [\operatorname{GCM}
\{
( G_{n}^{\star}(ck), V_{n}^{\star}(ck) ),
k\in
\mathcal{I}_{n}\}]
(G_{n}^{\star}(ci)).
\end{eqnarray}
For simplicity of notation, in the remainder of this section, we will
often write an integer interval
as a usual interval with two integer endpoints. This will, however, not
cause confusion since the interpretation
of the interval will be immediate from the context.

The following theorem gives the limit behavior of $\mathbb{X}_n$.
%
%th3.7 #&#
\begin{theorem}[(Weak convergence of $\mathbb{X}_{n}$)]
\label{thmgamma-033Xn-Ynfinite-dimensional-convergence}
Under the same assumptions as in Theorem \ref{thmgamma-033-1Xn-Yn},
for each nonnegative integer $N$, we have
\[
\{\mathbb{X}_{n}(ci), i\in[-N,N] \}
\stackrel{d}{\rightarrow}
\{\mathbb{X}_{c,\alpha,\beta} (ci), i\in[-N,N]\}.
\]
It follows that $n^{1/3}(\hat{F}(t_l) - F(t_l)) \stackrel
{d}{\rightarrow} \mathbb{X}_{c,\alpha,\beta}(0)$.
% Furthermore, we have
% \begin{eqnarray*}
% (\{\mathbb{X}_{n}(c\lfloor h/c \rfloor), h\in I_{h} \},
% \{\mathbb{Y}_{n}(c\lfloor h/c \rfloor), h\in I_{h} \})\\
% &\rightsquigarrow&
%(\{\mathbb{X}(c\lfloor h/c \rfloor), h\in I_{h} \},
%in the product space
%$L^{p}I_{h}\times L^{p}I_{h}$
%for each $p\geq1$,
% where
%$I_{h}=[-cC_{1}, cC_{2}]$
%with positive integers $C_{1}$ and $C_{2}$.
\end{theorem}
%
%re3.8 #&#
\begin{rem}
It is interesting to note the change in the limiting behavior of the
NPMLE with varying $\gamma$. As noted previously, for $\gamma\in(0,
1/3)$, the grid is sparse enough so that the naive average responses at
each inspection time, which provide
empirical estimates of $F$ at those corresponding inspection times, are
automatically ordered (and therefore the solution to the isotonic
regression problem) and there is no ``strength borrowed'' from nearby
inspection times. Consequently, a Gaussian limit is obtained. For
$\gamma\geq1/3$, the grid points are ``close enough,'' so that the
naive pointwise averages are no longer the best estimates of $F$. In
fact, owing to the closeness of successive grid-points, the naive
averages are no longer ordered, and the PAV pool adjacent violators
algorithm (PAVA) leads to a nontrivial solution for the NPMLE which is
a highly nonlinear functional of the data, putting us in the setting of
nonregular asymptotics. It turns out that for $\gamma\geq1/3$, the
order of the local neighborhoods of $t_l$ that determine the value of
$\hat{F}(t_l)$ is $n^{-1/3}$. When $\gamma= 1/3$, the resolution of
the grid matches the order of the local neighborhoods, leading in the
limit to a process in discrete-time that depends on $c$. When $\gamma>
1/3$, the number of grid-points in an $n^{-1/3}$ neighborhood of $t_l$
goes to infinity. This eventually washes out the dependence on $c$ and
also produces, in the limit, a~process in continuous time.
\end{rem}

For the rest of this section, we refer to the process $\mathbb
{X}_{c,\alpha,\beta}$ simply as $\mathbb{X}$ and the process
$\mathcal{P}_{c,\alpha,\beta}$ as $\mathcal{P}_c$.
\begin{pf*}{Proof--sketch of Theorem \ref
{thmgamma-033Xn-Ynfinite-dimensional-convergence}} The key steps of
the proof are as follows.
%First, consider the finite dimensional weak convergence.
%In fact, we again only need to follow the lines of the proof
%of Theorem \ref{thmgamma033Xn}.
%Note that the corresponding two claims and the first two facts can be
%established in the same way as those in the proof of
%Theorem \ref{thmgamma033Xn}.
%The third fact involves a little extension, which is established
%in Lemma \ref{Lemmagamma033finitedimensionalconvergencefact3}.
Take an integer $M>N$. Then, the following two claims hold.
\begin{Claim}\label{claim1}
There exist (integer-valued) random variables $L_{n} <-M$
and $U_{n} >M$ which are $O_{P}(1)$ and satisfy
\begin{eqnarray*}
&& \operatorname{GCM}
\{( G_{n}^{\star}(ck), V_{n}^{\star}(ck) ), k\in[L_{n},
U_{n}]
\}\\
&&\qquad= \operatorname{GCM} \{( G_{n}^{\star}(ck), V_{n}^{\star}(ck) ), k\in
\mathbb{Z} \}| [G_{n}^{\star}(cL_{n}), G_{n}^{\star}(cU_{n})].
\end{eqnarray*}
\end{Claim}
\begin{Claim}\label{claim2}
There also exist (integer-valued) random variables $L
<-M$ and $U >M$ such that
$L, U$ are $O_{P}(1)$ and that
\[
\operatorname{GCM}
\{\mathcal{P}_{c}(k), k\in[L, U]
\}= \operatorname{GCM} \{\mathcal{P}_{c}(k), k\in\mathbb{Z} \}| [cL,
cU].
\]
For the proofs of these claims, see \citet{Tang2011}. We next need a
key approximation lemma, which is a simple extension of Lemma~4.2 in
\citet{Rao1969}.
\end{Claim}
%For its proof, we only need to consider Levy-Prokhorov metric
%instead of Levy metric.
%
%le3.9 #&#
\begin{lem}\label{LemmaTruncation}
Suppose that for each $\varepsilon>0$,
$\{W_{n\varepsilon}\}$, $\{W_{n}\}$ and $\{W_{\varepsilon}\}$
are sequences of random vectors, $W$ is a random vector and that:
\begin{longlist}[(3)]
\item[(1)] $\lim_{\varepsilon\rightarrow0}
\overline{\lim}_{n\rightarrow\infty}
\mathbb{P}( W_{n\varepsilon} \not= W_{n} )= 0$,
\item[(2)] $\lim_{\varepsilon\rightarrow0}
\mathbb{P}( W_{\varepsilon} \not= W )= 0$,
\item[(3)] $W_{n\varepsilon} \stackrel{d}{\rightarrow} W_{\varepsilon}$,
as $n\rightarrow\infty$ for each $\varepsilon>0$.
\end{longlist}
Then $W_{n}\stackrel{d}{\rightarrow} W$, as $n\rightarrow\infty$.\vadjust{\goodbreak}
\end{lem}

From Claims \ref{claim1} and \ref{claim2}, for every (small) $\varepsilon> 0$,
there exists an integer $M_{\varepsilon}$ large enough
such that
\[
P( M_{\varepsilon} > \max
\{|L_{n}|,U_{n}, |L|,U\})
>1-\varepsilon.
\]
Denote, for $i\in[-N, N]$,
\begin{eqnarray*}
\mathbb{X}_{n}^{M_{\varepsilon}}(ci)
&=&
\operatorname{LS} \bigl[\operatorname{GCM}
\{( G_{n}^{\star}(ck), V_{n}^{\star}(ck) ), k\in
[\pm M_{\varepsilon}]\}\bigr]
(G_{n}^{\star}(ci)), \\
\mathbb{X}^{M_{\varepsilon}}(ci) &=&
\operatorname{LS} \bigl[\operatorname{GCM}
\{\mathcal{P}_{c}(k), k\in[\pm M_{\varepsilon}] \}\bigr]  (ci).
\end{eqnarray*}
Denote $[\pm N]=[-N, N]$ and
\begin{eqnarray*}
A_{n} &=& \bigl\{\{
\mathbb{X}_{n}^{M_{\varepsilon}}(ci), i \in[\pm N] \}
\not=
\{
\mathbb{X}_{n}(ci), i \in[\pm N] \}\bigr\}, \\
A &=& \bigl\{\{
\mathbb{X}^{M_{\varepsilon}}(ci), i \in[\pm N] \}
\not=
\{
\mathbb{X}(ci), i \in[\pm N] \}\bigr\}.
\end{eqnarray*}
Then, the following three facts hold:
\begin{Fact}\label{fact1} $\lim_{\varepsilon\rightarrow0}
\overline{\lim}_{n\rightarrow\infty}
\mathbb{P}
( A_{n} )= 0$.
\end{Fact}
\begin{Fact}\label{fact2}
$\lim_{\varepsilon\rightarrow0}
\mathbb{P}( A )= 0$.
\end{Fact}
\begin{Fact}\label{fact3}
$\{ \mathbb{X}_{n}^{M_{\varepsilon}}(ci), i \in[\pm N] \}\stackrel
{d}{\rightarrow} \{ \mathbb{X}^{M_{\varepsilon}}(ci), i \in[\pm N] \}$,
as $n\rightarrow\infty$ for each $\varepsilon>0$.
\end{Fact}

Facts \ref{fact1} and \ref{fact2} follow since $A_{n}$ and $A$
are subsets of
$\{ M_{\varepsilon} \leq\max\{|L_{n}|, U_{n},\allowbreak |L|, U \}\}$,
whose probability is less than $\varepsilon$,
Facts \ref{fact1} and \ref{fact2} hold. Fact \ref{fact3} is proved in \citet{Tang2011}. A direct
application of Lemma \ref{LemmaTruncation} then leads to the weak
convergence that we
sought to prove.
\end{pf*}
%
%re3.10 #&#
\begin{rem}
The proofs of Claims \ref{claim1} and \ref{claim2} consist of technically important
localization arguments.
Claim \ref{claim1} ensures that eventually, with arbitrarily high pre-specified
probability, the restriction of the greatest convex minorant of the
process $(G_n^{\star}, V_n^{\star})$ (which is involved in the
construction of~$\mathbb{X}_n$) to a bounded domain can be made equal
to the greatest convex minorant of the restriction of $(G_n^{\star},
V_n^{\star})$ to that domain, provided the domain is chosen
appropriately large, depending on the pre-specified probability. It can
be proved by using techniques similar to those in Section 6 of \citet{Kim1990}.
Claim \ref{claim2} ensures that an analogous phenomenon holds for the greatest
convex minorant of the process $\mathcal{P}_{c}$,
which is involved in the construction of $\mathbb{X}$.
These equalities then translate to the left-derivatives
of the GCMs involved, and the proof is completed
by invoking a continuous mapping theorem for the GCMs of
the restriction of $(G_n^{\star}, V_n^{\star})$ to
bounded domains, along with Claims \ref{claim1} and \ref{claim2},
which enable the use of the approximation lemma adapted from \citet{Rao1969}.

The basic strategy of the above proof has been invoked time and again
in the literature on monotone function estimation. \citet{Rao1969}
employed this technique to determine the limit distribution of the
Grenander estimator at a point, and \citet{Brunk1970} for studying
monotone regression. \citet{Leurgans1982} extended these techniques to
more general settings which cover weakly dependent data while \citet
{Anevski2006} provided a~comprehensive and unified treatment of
asymptotic inference under order restrictions, applicable to
independent as well as short and long range dependent data. This
technique was also used in \citet{Banerjee2007d} to study the
asymptotic distributions of a very general class of monotone response
models. It ought to be possible to bring the general techniques of
\citet{Anevski2006} to bear upon the boundary case, but we have not
investigated that option; our proof-strategy is most closely aligned
with the proof of Theorem~2.1 in \citet{Banerjee2007d}.
\end{rem}

%s3.4 #&#
\subsection{A brief discussion of the boundary phenomenon}
We refer to the behavior of the NPMLE for $\gamma= 1/3$ as the
\textit{boundary phenomenon}. As indicated in the
\hyperref[intro]{Introduction}, the asymptotic distribution for $\gamma=
1/3$ is \textit{different} from both the Gaussian (which comes into
play for $\gamma< 1/3$) and the Chernoff (which arises for $\gamma>
1/3$). This boundary distribution, which depends on the scale
parameter, $c$, can be viewed as an intermediate between the Gaussian
and Chernoff, and its degree of proximity to one or the other is
dictated by $c$ as we demonstrate in the following section. More
importantly, this transition from one distribution to another via the
boundary one, has important ramifications for inference in our
grid-based problem as also demonstrated in the next section.

The closest result to our boundary phenomenon in the literature appears
in the work of \citet{Zhang2001} who study the asymptotics of isotonic
estimation of a decreasing density with histogram-type data. Thus, the
domain of the density is split into a number of \textit{pre-specified
bins}, and the statistician knows the number of i.i.d. observations
from the density that fall into each bin (with a total of $n$ such
observations). The rate at which the number of bins increases relative
to $n$ then drives the asymptotics of the NPMLE of the density within
the class of decreasing piecewise linear
densities, with a distribution similar to $\mathbb{X}(0)$ appearing
when this number increases at rate $n^{1/3}$. However, \textit{unlike
us}, \citet{Zhang2001} do not establish any connections among the
different limiting regimes; neither do they offer a prescription for
inference when the rate of growth of the bins is unknown as is usually
the case in practice.

It is worthwhile contrasting our boundary phenomenon with those
observed by some other authors. \citet{Anevski2006} discover a
``boundary effect'' in their Theorems 5 and 6.1 when dealing with an
isotonized version of a kernel estimate (see Section 3.3 of their
paper). In the setting of i.i.d. data, when the smoothing bandwidth is
chosen to be of order $n^{-1/3}$, the asymptotics of the isotonized
kernel estimator are given by the minimizer of \textit{a Gaussian
process} (\textit{depending on the kernel}) \textit{with
continuous sample paths plus a quadratic drift}, whereas\vspace*{1pt} for bandwidths
of larger orders than $n^{-1/3}$ normal distributions obtain. A similar
phenomenon, in the setting of monotone density estimation, was observed
by \citet {Vaart2003} in their Theorem 2.2 for an isotonized
kernel estimate of a decreasing density while using an $n^{-1/3}$ order
bandwidth. Note that these boundary effects are quite different from
our boundary phenomenon. \textit{In Anevski and Hossjer's setting}, for
example, the underlying regression model is observed on the grid
$\{i/n\}$, with one response per grid-point. Kernel estimation with an
$n^{-1/3}$ bandwidth smooths the responses over time-neighborhoods of
order $n^{-1/3}$ producing a continuous estimator which is then
subjected to isotonization. This leads to a limit that is characterized
in terms of a process in \textit{continuous time}. \textit{In our
setting}, our data are not necessarily observed on an $\{i/n\}$ grid;
our grids can be much sparser and for the case $\gamma= 1/3$, multiple
responses are available at each grid-point. The NPMLE isotonizes the
$\bar {Z}_i$'s; thus, isotonization is preceded by averaging the
multiple responses at each time cross-section, but \textit{there is no
averaging of responses across time}, in sharp contrast to Anevski and
Hossjer's setting. This, in conjunction with the already noted fact at
the beginning of this subsection that the grid-resolution when $\gamma=
1/3$ has the same order as the localization involved in constructing
the process $\mathbb{X}_n$, leads in our case to a limit distribution
for the NPMLE that is characterized as a functional of a process in
\textit{discrete time}.

%%\newpage
%
%%\newpage

%s4 #&#
\section{Adaptive inference for $F$ at a point}\label{sec4}

In this section, we develop a~procedure for constructing asymptotic
confidence intervals for $F(t_{l})$
which does not require knowing or estimating the underlying grid
resolution controlled by the parameters $\gamma$ and $c$.
This provides massive advantage from an inferential perspective because
the parameter $\gamma$ critically drives the limit
distribution of the NPMLE and mis-specification of $\gamma$ may result
in asymptotically incorrect confidence sets, either due to
the use of the wrong limit distribution or due to an incorrect
convergence rate, or both.

To this end, we first investigate the relationships among the three
different asymptotic
limits for $\hat{F}(t_{l})$
that were derived in the previous section, for different values of
$\gamma$. In what follows, we denote $\mathbb{X}_{c,\alpha,\beta
}(0)$ by $\mathcal{S}_c$, suppressing the dependence
on $\alpha,\beta$ for notational convenience. The use of the letter
$\mathcal{S}$ is to emphasize the characterization of this random
variable as the
slope of a stochastic process.

Our first result relates the distribution of $\mathcal{S}_c$ to the Gaussian.
%
%th4.1 #&#
\begin{theorem}\label{thmrelationshipfrom-one-third-to-lessWald}
As $c\rightarrow\infty$,
% \sqrt{c}\mathcal{S}_c
% \stackrel{d}{\rightarrow} \alpha Z,
$
\sqrt{c}\mathcal{S}_c
\stackrel{d}{\rightarrow} \alpha Z,
$
where $Z$ follows the standard normal distribution.
\end{theorem}

Our next result investigates the case where $c$ goes to 0.
\begin{theorem}\label{thmrelationshipfrom-one-third-to-moreWald}
As $c\rightarrow0$,
% \mathcal{S}_{c}
% \stackrel{d}{\rightarrow}
% g_{\alpha,\beta}(0)
% \stackrel{d}{=}2(\alpha^{2}\beta)^{1/3} \mathcal{Z}.
$
\mathcal{S}_{c}
\stackrel{d}{\rightarrow}
g_{\alpha,\beta}(0)
\stackrel{d}{=}2(\alpha^{2}\beta)^{1/3} \mathcal{Z}.
$
\end{theorem}
%
%re4.3 #&#
\begin{rem}
%Comparing with Theorem \ref{thmgammalarger033distributionWald}.
Theorem \ref{thmrelationshipfrom-one-third-to-moreWald} is somewhat
easier to visualize heuristically, compared to
Theorem \ref{thmrelationshipfrom-one-third-to-lessWald}.
Recall that $\mathcal{S}_c$ is the left-slope of
the GCM of the process $\mathcal{P}_c$ at the point 0,
the process itself being defined on the grid $c\mathbb{Z}$.
As $c$ goes to 0,
the grid becomes finer, and
the process $\mathcal{P}_c$ is eventually substituted
by its limiting version,
namely $X_{\alpha,\beta}$.
Thus, in the limit, $\mathcal{S}_c$ becomes $g_{\alpha,\beta}(0)$,
the left-slope of the GCM of $X_{\alpha,\beta}$ at 0.
The representation of this limit in terms of $\mathcal{Z}$
was established in
Remark
\ref{remgamma-greaterThan-oneThirdspecial}
following Theorem \ref{thmgamma-033-1Xn-Yn}.
\end{rem}

%f1 #&#
\begin{figure}

\includegraphics{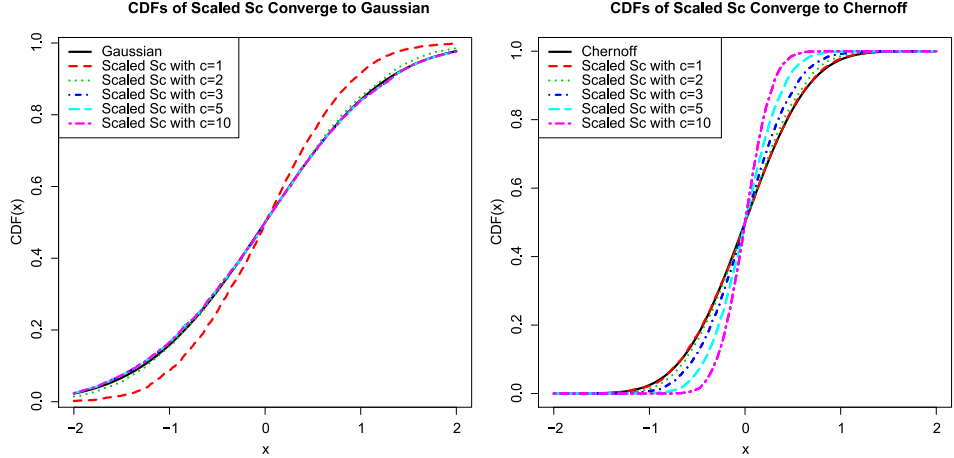}

\caption{The left and right panels show that a sequence of empirical
CDFs of
the properly scaled $\mathcal{S}_{c}$ converge to the standard
Gaussian and Chernoff distributions,
respectively.
In the left panel,
the empirical CDFs with $c \geq3$ almost coincide with the standard
Gaussian distribution.}
\label{fgScConvergeGassianChernoff}
\end{figure}

The results of
Theorems
\ref{thmrelationshipfrom-one-third-to-lessWald}
and
\ref{thmrelationshipfrom-one-third-to-moreWald} are illustrated next.
Suppose the time interval $[a,b]$ is [0,2], $x_{0}=1$ and that $F$ and
$G$ are
both the uniform distribution on $[0,2]$.
Under these settings, the values of $\alpha$ and $\beta$ are
$\sqrt{2}/4$ and $1/4$, respectively.
We generate i.i.d. random samples of $\mathcal{S}_{c}$
with $c$ being 1, 2, 3, 5 and 10 and the common sample size being 5000.
The left panel of Figure \ref{fgScConvergeGassianChernoff}
compares the empirical cumulative distribution functions (CDF)
of $\sqrt{c}\mathcal{S}_{c}/\alpha$
and the standard Gaussian distribution $N(0,1)$.
It shows clearly that the empirical CDFs move closer to the Gaussian
distribution
with increasing~$c$ and that the empirical CDF of $\sqrt{c}\mathcal
{S}_{c}/\alpha$
with $c$ equal to 3 has already provided a decent approximation to $N(0,1)$.
On the other hand, the right panel of Figure \ref
{fgScConvergeGassianChernoff}
compares the empirical CDFs
of $(1/2)(\alpha^{2}\beta)^{-1/3} \mathcal{S}_{c}$
and the standard Chernoff distribution $\mathcal{Z}$.
Again, the empirical CDFs approach that of $\mathcal{Z}$ with
diminishing $c$, with $c=1$
providing a close approximation for $\mathcal{Z}$. Note that, while
the convergence in this setting is relatively quick in the sense
that the limiting phenomena manifest themselves at moderate values of
$c$ (i.e., neither too large, nor too small), this may not necessarily
be the
case for other combinations of $(\alpha,\beta)$, and more extreme
values may be required for good enough approximations.\vspace*{8pt}

%So, with the specific $\alpha$ and $\beta$ used in the simulation,
%$\mathcal{S}_{c}$ converges to Gaussian
%without requiring c to be very large;
%at the same time,
%$\mathcal{S}_{c}$ converges to Chernoff
%without requiring c to be very small.
%This means that
%$\mathcal{S}_{c}$
%usually could provide good approximation for
%both Gaussian and Chernoff
%for the case in the simulation.
%Obviously, it is because that
%the values of $\alpha$ and $\beta$ used in the simulation
%come from a simple underlying model.

\textit{The adaptive inference scheme}: We are now in a position to
propose our inference scheme.
We focus on the so-called ``Wald-type'' intervals for $F(t_{l})$, that is,
intervals of the form $\hat{F}(t_{l})$ plus and minus terms
depending on the sample size and
the large sample distribution of the estimator.
Let $c_0$ and~$\gamma_0$ denote the \textit{true unknown values} of $c$
and~$\gamma$ in the current status model. With
$K = K_n$ being the number of grid-points, we have the relation
\[
K_n = \lfloor{(b-a)/(c_0 n^{-\gamma_0})}\rfloor .
\]
Now \textit{pretend} that the true $\gamma$ is exactly equal to $1/3$.
Calculate a surrogate~$c$, say $\hat{c}$, via the relation
\[
\lfloor{(b-a)/(\hat{c}n^{-1/3})}\rfloor= K_n  .
\]
%
%Since the `floors' in the above relation
%between $c$ and $\hat{c}$ are asymptotically ignorable,
%we theoretically specify $\hat{c}$ to be $cn^{1/3-\gamma}$.
%Thus, the calculated parameter $\hat{c}$ actually depends on $n$,
%and goes to $\infty$ and 0
%for $\gamma\in(0,1/3)$ and $\gamma\in(1/3,1)$, respectively.
%Henceforth, we denote $\tilde{c}$ by $\hat{c}$.
Some algebra shows that
\[
\hat{c} = \hat{c}_n = c n^{1/3 - \gamma_0} + O(n^{1/3 - 2 \gamma
_0}) = c n^{1/3 - \gamma_0}\bigl(1 + O(n^{-\gamma_0})\bigr)  .
\]
Thus, the calculated parameter $\hat{c}$ actually depends on $n$, and
goes to $\infty$ and 0 for $\gamma_0 \in(0,1/3)$ and $\gamma_0 \in
(1/3,1]$, respectively.

We propose\vspace*{1pt} to use the distribution of $\mathcal{S}_{\hat{c}}$
as an approximation to the distribution of $n^{1/3}(\hat{F}(t_{l}) -
F(t_{l}))$.
%as is suggested by our `boundary asymptotics' results (REF),
Thus, an adaptive approximate
$1-\eta$
confidence interval for $F(t_{l})$
is given by
%
%e4.1 #&#
\begin{equation}\label{CIWaldadaptive}
\bigl[\hat{F}(t_{l}) - n^{-1/3} q(\mathcal{S}_{\hat{c}}, 1 - \eta/2),
\hat{F}(t_{l}) - n^{-1/3} q\bigl(\mathcal{S}_{\hat{c}}, (\eta/2)\bigr)\bigr],
\end{equation}
where $\eta>0$ and $q(X,p)$ stands for the lower $p$th quantile of a
random variable $X$
with $p\in(0,1)$.\vspace*{8pt}

\textit{Asymptotic validity of the proposed inference scheme}: The above
adaptive confidence interval
provides the correct asymptotic calibration, \textit{irrespective of the
true value of $\gamma$}.
If $\gamma_0$ happens to be $1/3$, then, of course, the adaptive
confidence interval
is constructed with the correct asymptotic result. If not, consider
first the case that $\gamma_0 \in(1/3, 1]$.
%In the case that the true $\gamma> 1/3$, this is easy to see.
If \textit{we knew} that $\gamma_0 \in(1/3,1]$,
then, by result (\ref{resultgamma-033-1Xn0})
and the symmetry of $g_{\alpha,\beta}(0)$,
the true confidence interval would be
%
%e4.2 #&#
\begin{equation}\label{CIWaldgamma-033-1}
\bigl[\hat{F}(t_{l}) \pm
n^{-1/3} q\bigl(g_{\alpha,\beta}(0), (1-\eta/2)\bigr)\bigr].
\end{equation}
Now recall that $\hat{c}$ goes to 0 since $\gamma_0\in(1/3,1]$.
Thus, by Theorem \ref{thmrelationshipfrom-one-third-to-moreWald},
the quantile sequence $q(\mathcal{S}_{\hat{c}}, p)$ converges to
$q(g_{\alpha,\beta}(0), p)$,
owing to the fact that~$g_{\alpha,\beta}(0)$ is a continuous random variable.
So, the adaptive confidence interval~(\ref{CIWaldadaptive})
converges to the true one (\ref{CIWaldgamma-033-1})
obtained when $\gamma_0$ is in $(1/3,1]$.\vadjust{\goodbreak}

That the adaptive procedure also works when $\gamma_0 \in(0,1/3)$
will be shown by using Theorem \ref
{thmrelationshipfrom-one-third-to-lessWald}.
Again, \textit{suppose we know} the value of $\gamma_0$. Then,
from result (\ref{equcase1hat-Flimiting-distribution2nd-formation})
and the symmetry of the standard normal random variable $Z$,
the confidence interval is given by
%
%e4.3 #&#
\begin{equation}\label{CIWaldgamma-0-033}
\bigl[\hat{F}(t_{l}) \pm n^{-(1-\gamma_0)/2}\alpha c^{-1/2} q\bigl(Z,
(1-\eta/2)\bigr) \bigr].
\end{equation}
To show that the adaptive procedure is, again, asymptotically correct,
it suffices to show that
for every $p\in(0,1)$,
as $n \rightarrow\infty$,
\[
\frac{n^{-1/3} q(\mathcal{S}_{\hat{c}}, p)}
{n^{-(1-\gamma_0)/2}\alpha c^{-1/2} q(Z, p)}
=
\frac{n^{-1/3}c^{1/2}}
{n^{-(1-\gamma_0)/2}\hat{c}^{1/2}}
\cdot
\frac{\hat{c}^{1/2}q(\mathcal{S}_{\hat{c}}, p)}
{\alpha q(Z, p)}
=
I \cdot \mathit{II}
\rightarrow1.
\]
Recall that $\hat{c}$ goes to $\infty$ since $\gamma_0\in(0,1/3)$.
By Theorem \ref{thmrelationshipfrom-one-third-to-lessWald},
we have $\mathit{II}\rightarrow1$ as $n \rightarrow\infty$.
On the other hand, we can see $I$ simplifies to $(1+ O(n^{-\gamma
_0}))^{-1/2}$ and therefore goes to 1.
Thus, the adaptive confidence interval (\ref{CIWaldadaptive})
also converges to
the true one (\ref{CIWaldgamma-0-033})
obtained when $\gamma_0$ is known to be in $(0, 1/3)$.

Thus, our procedure \textit{adjusts automatically} to the inherent rate of
growth of the number of distinct observation times and that is an
extremely desirable property.

We next articulate some practical issues with the adaptive procedure.
First, note that $\mathcal{S}_{\hat{c}} = \mathbb{X}_{\hat
{c},\alpha,\beta}(0)$,
and in practice $\alpha$ and $\beta$ are unknown, and therefore need
to be estimated consistently. We provide simple methods for consistent
estimation of these two parameters in the next section. Second, the
random variable $\mathbb{X}_{\hat c,\alpha,\beta}(0)$ does not
appear to
admit a natural scaling in terms of some canonical random variable: in
other words, it cannot be represented as $C(c,\alpha,\beta) J$ where
$C$ is an explicit function of $c,\alpha,\beta$ and~$J$ is some fixed
well-characterized random variable.
Thus, the quantiles of~$\mathbb{X}_{\hat{c},\hat{\alpha},\hat
{\beta}}$ (where~$\hat{\alpha}$ and\vspace*{1pt} $\hat{\beta}$ are consistent
estimates for the corresponding parameters) need to be calculated by
generating many sample paths
from the parent process $\mathcal{P}_{\hat{c},\hat{\alpha},\hat
{\beta}}$ and computing the left slope of the convex minorant of each
such path at 0.
This is, however, not a terribly major issue in these days of fast
computing, and, in our opinion,
the mileage obtained in terms of adaptivity more than compensates for
the lack of scaling. Finally, one may wonder if resampling the NPMLE
would allow adaptation with respect to~$\gamma$. The problem, however,
lies in the fact that while the usual~$n$ out of $n$ bootstrap works for
the NPMLE when $\gamma\in(0,1/3)$, it fails under the nonstandard
asymptotic regimes that operate for
$\gamma\in[1/3,1]$, as is clear from the work of \citet{Abrevaya2005},
\citet{Kosorok2008a} and \citet{Sen2010}. Since $\gamma$ is unknown,
it is impossible to decide
whether to use the standard~$n$ out of $n$ bootstrap. One could argue
that the~$m$ out of~$n$ bootstrap or subsampling
will work irrespective of the value of~$\gamma$, but the problem that
arises here is that
these procedures require knowledge of the convergence rate and this is
unknown as it depends on the true value of~$\gamma$.

%%\newpage

%s5 #&#
\section{A practical procedure and simulations}\label{sec5}

In this section, we provide a~practical version of the adaptive
procedure introduced in Section \ref{sec4} to construct Wald-type confidence
intervals for $F(t_{l})$ and assess their performance through
simulation studies.
%For the rest of this section, the
%process $\mathbb{X}_{c,\kappa_1,\kappa_2}$ will be abbreviated to $
The true values of $c$ and $\gamma$ are denoted by $c_0$ and $\gamma
_0$. The process
$\mathcal{P}_{c,\alpha,\beta}$ is again abbreviated to $\mathcal{P}_c$.

Recall\vspace*{1pt} that in the adaptive procedure,
we always specify $\gamma=1/3$
and compute a surrogate for $c_0$, namely $\hat c$, as a solution of
the equation
$K = \lfloor{(b-a)/\hat{c}n^{-1/3}\rfloor}$,
where $K$ is the number of grid points.
To construct a level $1-2\eta$ confidence interval for $F(t_{l})$ for
a small positive $\eta$,
quantiles of $\mathcal{S}_{\hat{c}}$ are needed.
Since $\mathcal{S}_c = \operatorname{LS}[\operatorname{GCM} \{\mathcal{P}_c(k), k\in\mathbb{Z}
\}] (0)$
($c$ is genetically used),
we approximate $\mathcal{S}_c$ with
\[
\mathbb{X}_{c,K_{a}}(0)= \operatorname{LS}\bigl[\operatorname{GCM}
\{\mathcal{P}_{c}(k), k\in[-K_{a}-1, K_{a}] \}\bigr] (0)
\]
for some large $K_{a}\in\mathbb{N}$.
Further, since
\[
\mathbb{X}_{c,K_{a}}(0)
=\operatorname{LS}\bigl[\operatorname{GCM}
\bigl\{\bigl(\mathcal{P}_{1,c}(k)/c, \mathcal{P}_{2,c}(k)/c\bigr), k\in
[-K_{a}-1, K_{a}] \bigr\}\bigr] (0),
\]
where $\mathcal{P}_{1,c}(k)/c=k$ and
$\mathcal{P}_{2,c}(k)/c = \alpha W(ck)/c + \beta ck(1+k)$,
we get that $\mathbb{X}_{c,K_{a}}(0)$ is the isotonic regression at $k=0$
of the data
\begin{eqnarray*}
&& \bigl\{\bigl(k, \mathcal{P}_{2,c}(k)/c-\mathcal{P}_{2,c}(k-1)/c\bigr), k\in
[-K_{a}, K_{a}] \bigr\}\\
&&\qquad= \bigl\{\bigl(k,\alpha Z_{k}/\sqrt{c} + 2\beta ck\bigr), k\in[-K_{a}, K_{a}]
\bigr\},
\end{eqnarray*}
where\vspace*{1pt} $\{Z_{k}\}_{k=-K_{a}}^{K_{a}}$ are i.i.d. from
$N(0,1)$, $\alpha= \sqrt{F(x_{0})(1-F(x_{0}))/g(x_{0})}$ and $\beta=
f(x_{0})/2$. To make this adaptive procedure practical, we next
consider the estimation of $\alpha$ and $\beta$, or equivalently, the
estimation of $F(x_{0})$, $g(x_{0})$ and $f(x_{0})$.

First, we consider the estimation of $F(x_{0})$ and $g(x_{0})$.
Although
$F(x_{0})$ can be consistently estimated by $\hat F(t_{l})$,
in our simulations we estimate~$F(x_{0})$ by
$\rho\hat F(t_{l}) + (1-\rho)\hat F(t_{r})$
with $\rho=(x_{0}-t_{l})/(t_{r}-t_{l})\in[0,1)$.
To estimate~$g(x_{0})$,
we use
the following estimating equation:
$
(N_{l-j^{\star}+1} + \cdots+ N_{r+j^{\star}})/n
=\break g(x_{0})(t_{r+j^{\star}}-t_{l-j^{\star}})
$,
where $j^{\star}$ is defined below in the estimation of $f(x_{0})$.
Since the design density $g$ is assumed to be continuous in a
neighborhood of $x_{0}$,
and the interval $[t_{l-j^{\star}}, t_{r+j^{\star}}]$ is shrinking to $x_{0}$,
it is reasonable to approximate $g$ over
the interval $[t_{l-j^{\star}}, t_{r+j^{\star}}]$ with a constant function.
Thus, from the above estimating equation,
one simple but consistent estimator of~$g(x_{0})$ is given by
$
\hat g(x_{0}) =(N_{l-j^{\star}+1} + \cdots+ N_{r+j^{\star}})
/[n(t_{r+j^{\star}}-t_{l-j^{\star}})]
$.

Next, we consider the estimation of $f(x_{0})$. To this end, we
estimate $f(t_{l})$ using a local linear approximation:
identify a small interval around $t_{l}$, and then approximate $F$ over this
interval by a line, whose slope gives the estimator of $f(t_{l})$.
We determine the interval by
the following several requirements.
First, the sample proportion $p_{n}$
in the interval should be larger than
the sample proportion at each grid point,
which is of order $n^{-\gamma}$ for $\gamma\in(0,1]$.
For example, setting $p_{n}$ be of order $1/\log n$ theoretically
ensures a sufficiently large interval.
Second, for simplicity, we make the interval symmetric around~$t_{l}$.
Third, in order to obtain\vadjust{\goodbreak}
a positive estimate [since $f(t_{l})$ is positive],
we symmetrically enlarge the interval satisfying the above two requirements
until the values of $\hat F$ at the two ends of the interval become different.
Thus, we first find $j^{\star}$, the smallest integer such that
$
\sum_{i=l-j^{\star}}^{l+j^{\star}}N_{i}/n \geq1/\log n
$.
Next,
we find $i^{\star}$, the smallest integer larger than $j^{\star}$
such that
$
\hat F(t_{l-i^{\star}}) < \hat F(t_{l+i^{\star}})
$ and employ a linear approximation over $[t_{l-i^{\star}},
t_{l+i^{\star}}]$.
More specifically,
we compute
\[
(\hat\beta_{0},\hat\beta_{1})
= \mathop{\arg\max}_{(\beta_{0}, \beta_{1})\in\mathbb{R}^{2}}
\Biggl\{\sum_{i=l-i^{\star}}^{l+i^{\star}}
\bigl(\hat F(t_{i}) - \beta_{0} - \beta_{1}t_{i}\bigr)^{2}N_{i} \Biggr\}
\]
and estimate $f(t_{l})$ [and $f(x_{0})$] by $\hat\beta_{1}$.
Once these nuisance parameters have been estimated,
the practical adaptive procedure can be implemented.

The above procedures provide consistent estimates of $g(x_0)$ and
$f(x_0)$ under the assumption of
a single derivative for $F$ and $G$ in a neighborhood of $x_0$,
irrespective of the value of $\gamma$ [since the estimates are
obtained by local polynomial fitting over a neighborhood of logarithmic
order (in $n$) around $x_0$ and such neighborhoods are guaranteed
to be asymptotically wider than $n^{-\gamma}$ for any $0 < \gamma\leq
1$]. Two points need to be noted. First, the $1/\log n$ threshold used to
determine $j^{\star}$ in the previous paragraph may need to be changed
to a multiple of $1/\log n$, depending on the sample size and the length
of the time interval. Second, the locally constant estimate of $g(x_0)$
discussed above could be replaced by a local linear (or quadratic)
estimate of
$g$, if the data strongly indicate that $G$ is changing sharply in a
neighborhood of~$x_0$.

To evaluate the finite sample performance of the practical adaptive
procedure, we also provide simulated confidence intervals
of an idealized (theoretical) adaptive procedure where the true values
of the
parameters $F(x_{0}), g(x_{0})$ and $f(x_{0})$ are used, but $\gamma$
is still practically assumed to be $1/3$,
and $c$ is taken as the previous~$\hat c$.
These confidence intervals can be considered as
the best Wald-type confidence intervals
based on the adaptive procedure.

%For comparison purpose,
%we also consider an ideal procedure
%where the true value of $\gamma$ is supposed to be known
%and the corresponding Wald-type asymptotic distribution
%is used to construct confidence intervals.
%In other words,
%for $\gamma\in(0,1/3)$ and $\gamma\in(1/3,1)$, we use the following
%two results:
% \sqrt{N_{l}}(\hat F(t_{l}) - F(t_{l}))
% \stackrel{d}{\rightarrow}
% [F(x_{0})(1-F(x_{0}))]^{1/2}N(0,1),
% n^{1/3}(\hat F(t_{l}) - F(t_{l}))
% \stackrel{d}{\rightarrow}
% (\frac{4f(x_{0})F(x_{0})(1-F(x_{0}))}{g(x_{0})})^{1/3}
%Other unknown parameters
%are estimated as before.

The simulation settings are as follows:
The sampling interval $[a,b]$ is $[0,1]$.
The design density $g$ is uniform on $[a,b]$.
The distribution of $T$ is
the uniform distribution over $[a,b]$ or
the exponential distribution with $\lambda=1$ or 2.
The anchor-point $x_{0}$ is 0.5.
The pair of grid-parameters $(\gamma, c)$ takes values
$(1/6, 1/6)$, $(1/4, 1/4)$, $(1/3, 1/2)$, $(1/2, 1)$,
$(2/3, 2)$ and $(3/4, 3)$.
The sample size $n$ ranges from 100 to 1000 by 100.
When generating the quantiles of~$\mathbb{X}_{\hat c}(0)$,
$K_{a}$ is set to be 300 and the corresponding iteration number 3000.
We are interested in constructing 95\% confidence intervals for $F(t_{l})$.
The iteration number for each simulation is 3000.

%f2 #&#
\begin{figure}%[!t]

\includegraphics{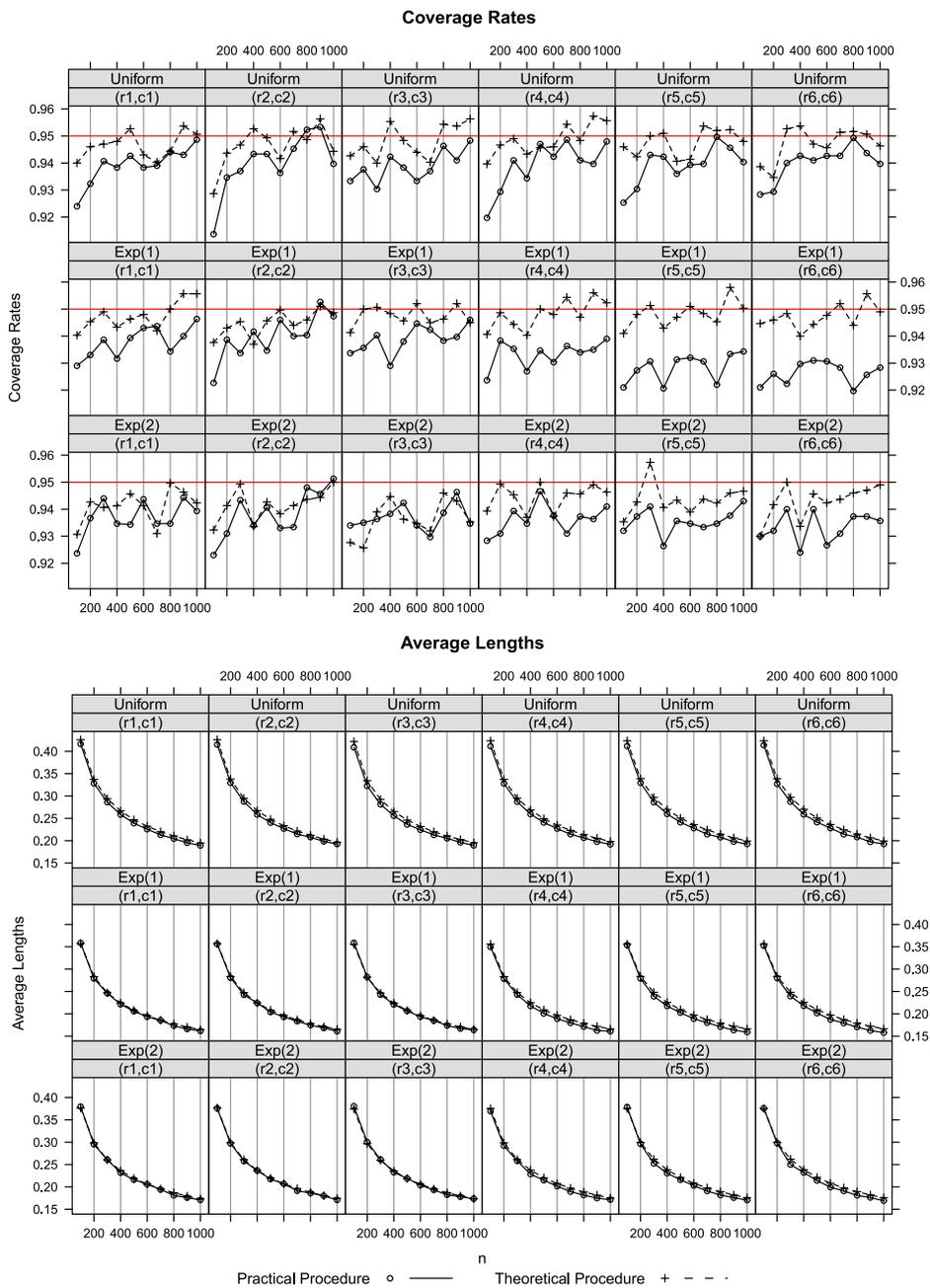}

\caption{A comparison of the coverage rates and average lengths
of the practical and theoretical procedures,
where $(ri, ci)$ for $i=1,\ldots,6$
are $(1/6, 1/6)$, $(1/4, 1/4)$, $(1/3, 1/2)$,
$(1/2, 1)$, $(2/3, 2)$ or $(3/4, 3)$,
respectively. The sample size $n$ varies from 100 to 1000 by 100.}
\label{fgsimuplotcr-al}
\end{figure}

%t1 #&#
\begin{table}[t!]
\caption{A comparison of the coverage rates and average lengths of the
practical procedure with those of the theoretical procedure, where
$U[0,1]$ and $\exp(\lambda)$ stand for the uniform distribution over
$[0,1]$, and the exponential distributions with the parameter
$\lambda$, and $n_{1}, n_{2}$ and $n_{3}$ are 100, 300 and 500,
respectively} \label{tablepractical-idea-proceduresCR-AL}
\vspace*{-3pt}
\begin{tabular*}{\tablewidth}{@{\extracolsep{\fill}}lccccccccc@{}}
\hline
\multicolumn{10}{@{}c@{}}{\textbf{Coverage rates}} \\
\hline
\textbf{CR(P)} & \multicolumn{3}{c}{$\bolds{U[0,1]}$}
& \multicolumn{3}{c}{$\bolds{\exp(1)}$} &
\multicolumn{3}{c@{}}{$\bolds{\exp(2)}$} \\[-4pt]
\multicolumn{1}{@{}c}{\hrulefill} & \multicolumn{3}{c}{\hrulefill}
& \multicolumn{3}{c}{\hrulefill} & \multicolumn{3}{c@{}}{\hrulefill}\\
$\bolds{(\gamma, c)}$ & $\bolds{n_{1}}$ & $\bolds{n_{2}}$
& $\bolds{n_{3}}$ & $\bolds{n_{1}}$ & $\bolds{n_{2}}$ &
$\bolds{n_{3}}$ & $\bolds{n_{1}}$ & $\bolds{n_{2}}$ & $\bolds{n_{3}}$ \\
\hline
$(1/6, 1/6)$ & 0.924 & 0.941 & 0.943 & 0.929 & 0.939 & 0.939 & 0.924 & 0.944 &
0.934 \\
$(1/4, 1/4)$ & 0.914 & 0.937 & 0.943 & 0.923 & 0.934 & 0.935 & 0.923 & 0.943 &
0.941 \\
$(1/3, 1/2)$ & 0.933 & 0.930 & 0.938 & 0.934 & 0.940 & 0.938 & 0.934 & 0.936 &
0.942 \\
$(1/2, 1)$ & 0.920 & 0.941 & 0.947 & 0.924 & 0.935 & 0.935 & 0.928 & 0.939 & 0.947
\\
$(2/3, 2)$ & 0.925 & 0.943 & 0.936 & 0.921 & 0.931 & 0.931 & 0.932 & 0.941 & 0.936
\\
$(3/4, 3)$ & 0.928 & 0.940 & 0.941 & 0.921 & 0.922 & 0.931 & 0.930 & 0.940 & 0.940\\
\end{tabular*}

\begin{tabular*}{\tablewidth}{@{\extracolsep{\fill}}lccccccccc@{}}
\hline
% after \  \hline or \cline{col1-col2} \cline{col3-col4} ...
\textbf{CR(T)} & \multicolumn{3}{c}{$\bolds{U[0,1]}$}
& \multicolumn{3}{c}{$\bolds{\exp
(1)}$} & \multicolumn{3}{c@{}}{$\bolds{\exp(2)}$} \\
\hline
%$(\gamma, c)$ & $n_{1}$ & $n_{2}$ & $n_{3}$ & $n_{1}$ & $n_{2}$ &
%$n_{3}$ & $n_{1}$ & $n_{2}$ & $n_{3}$ \\ \hline
$(1/6, 1/6)$ & 0.940 & 0.947 & 0.953 & 0.940 & 0.949 & 0.946 & 0.931 & 0.941 &
0.946 \\
$(1/4, 1/4)$ & 0.929 & 0.947 & 0.949 & 0.938 & 0.945 & 0.946 & 0.932 & 0.949 &
0.943 \\
$(1/3, 1/2)$ & 0.943 & 0.940 & 0.948 & 0.941 & 0.951 & 0.946 & 0.928 & 0.939 &
0.936 \\
$(1/2, 1)$ & 0.940 & 0.949 & 0.946 & 0.941 & 0.944 & 0.950 & 0.939 & 0.945 & 0.950
\\
$(2/3, 2)$ & 0.946 & 0.950 & 0.941 & 0.941 & 0.951 & 0.947 & 0.935 & 0.957 & 0.943
\\
$(3/4, 3)$ & 0.939 & 0.953 & 0.947 & 0.945 & 0.948 & 0.944 & 0.930 & 0.950 & 0.946
\end{tabular*}
\begin{tabular*}{\tablewidth}{@{\extracolsep{\fill}}lccccccccc@{}}
\hline
\multicolumn{10}{@{}c@{}}{\textbf{Average lengths}} \\
\hline
\textbf{AL(P)} & \multicolumn{3}{c}{$\bolds{U[0,1]}$}
& \multicolumn{3}{c}{$\bolds{\exp
(1)}$} & \multicolumn{3}{c@{}}{$\bolds{\exp(2)}$} \\[-4pt]
\multicolumn{1}{@{}c}{\hrulefill} & \multicolumn{3}{c}{\hrulefill}
& \multicolumn{3}{c}{\hrulefill} & \multicolumn{3}{c@{}}{\hrulefill}\\
$\bolds{(\gamma, c)}$ & $\bolds{n_{1}}$ & $\bolds{n_{2}}$
& $\bolds{n_{3}}$ & $\bolds{n_{1}}$ & $\bolds{n_{2}}$ &
$\bolds{n_{3}}$ & $\bolds{n_{1}}$ & $\bolds{n_{2}}$ & $\bolds{n_{3}}$ \\
\hline
$(1/6, 1/6)$ & 0.417 & 0.286 & 0.239 & 0.358 & 0.246 & 0.206 & 0.380 & 0.261 &
0.216 \\
$(1/4, 1/4)$ & 0.415 & 0.287 & 0.240 & 0.356 & 0.242 & 0.204 & 0.376 & 0.258 &
0.218 \\
$(1/3, 1/2)$ & 0.409 & 0.281 & 0.236 & 0.359 & 0.243 & 0.207 & 0.381 & 0.258 &
0.219 \\
$(1/2, 1)$ & 0.411 & 0.287 & 0.241 & 0.350 & 0.243 & 0.201 & 0.370 & 0.258 & 0.215
\\
$(2/3, 2)$ & 0.411 & 0.286 & 0.241 & 0.354 & 0.239 & 0.202 & 0.379 & 0.253 & 0.216
\\
$(3/4, 3)$ & 0.414 & 0.287 & 0.241 & 0.352 & 0.239 & 0.202 & 0.376 & 0.250 & 0.214
\end{tabular*}
\begin{tabular*}{\tablewidth}{@{\extracolsep{\fill}}lccccccccc@{}}
\hline
% after \  \hline or \cline{col1-col2} \cline{col3-col4} ...
\textbf{AL(T)} & \multicolumn{3}{c}{$\bolds{U[0,1]}$}
& \multicolumn{3}{c}{$\bolds{\exp
(1)}$} & \multicolumn{3}{c@{}}{$\bolds{\exp(2)}$} \\
\hline
%$(\gamma, c)$ & $n_{1}$ & $n_{2}$ & $n_{3}$ & $n_{1}$ & $n_{2}$ &
%$n_{3}$ & $n_{1}$ & $n_{2}$ & $n_{3}$ \\ \hline
$(1/6, 1/6)$ & 0.426 & 0.294 & 0.247 & 0.357 & 0.247 & 0.208 & 0.377 & 0.260 &
0.219 \\
$(1/4, 1/4)$ & 0.426 & 0.295 & 0.248 & 0.357 & 0.247 & 0.208 & 0.377 & 0.261 &
0.220 \\
$(1/3, 1/2)$ & 0.422 & 0.292 & 0.246 & 0.355 & 0.246 & 0.208 & 0.374 & 0.260 &
0.219 \\
$(1/2, 1)$ & 0.424 & 0.295 & 0.249 & 0.356 & 0.247 & 0.209 & 0.375 & 0.261 & 0.220
\\
$(2/3, 2)$ & 0.424 & 0.297 & 0.251 & 0.356 & 0.248 & 0.209 & 0.375 & 0.262 & 0.221
\\
$(3/4, 3)$ & 0.424 & 0.297 & 0.251 & 0.356 & 0.248 & 0.209 & 0.375 & 0.262 & 0.221
\\
\hline
\end{tabular*}\vspace*{-8pt}
\end{table}

Denote the simulated coverage rates and average lengths for the
practical procedure as CR(P) and AL(P)
and those for the theoretical procedure as CR(T) and AL(T).
Figure \ref{fgsimuplotcr-al}
contains the plots of
CR(P), CR(T), AL(P) and AL(T),
and Table \ref{tablepractical-idea-proceduresCR-AL}
contains the corresponding numerical values
for $n = 100, 300, 500$. The first panel of Figure \ref{fgsimuplotcr-al}
shows that both CR(T) and CR(P) are\vadjust{\goodbreak} usually close to
the nominal level 95\% from below
and that CR(T) is generally about 1\% better than CR(P).
This reflects the price of not knowing
the true values of the parameters
$F(x_{0})$, $g(x_{0})$ and $f(x_{0})$
in the practical procedure.
On the other hand,
the second panel of Figure \ref{fgsimuplotcr-al}
shows that the AL(P)s are usually slightly shorter than AL(T)s.
This indicates that the practical procedure is
slightly more aggressive.
As the sample size increases, the coverage rates usually approach the
nominal level,
and the average lengths also become shorter, as expected.

The patterns noted above show up in more extensive simulation studies,
not shown here owing to constraints of
space. Also, the adaptive procedure is seen to compete well with the
asymptotic approximations that one would use
for constructing CIs \textit{were} $\gamma$ \textit{known}.

We end this section by pointing out that while, for the simulations, we
knew the anchor-point $x_0$ ($t_l$ being the largest grid-point to the
left of or equal to $x_0$), and that we did make use of its value for
estimating $F(x_0)$ in our simulations, knowledge of $x_0$ is not
essential to the inference procedure. We could have just estimated
$F(x_0)$ by $\hat{F}(t_l)$ [rather than by a convex combination of
$\hat{F}(t_l)$ and
$\hat{F}(t_r)$ that depends upon $x_0$] consistently. This is a
critical observation, since in a real-life situation what we are
provided is current status data on a grid with particular grid points
of interest. There is no specification of $x_0$. To make inference on
the value of $F$ at such a grid-point, one can, conceptually, view
$x_0$ as being any point strictly in between the given point and the
grid-point immediately after, but its value is not required to
construct a confidence interval by the adaptive method. To reiterate,
the ``anchor-point,'' $x_0$ was introduced for developing our
theoretical results, but its value can be ignored for the
implementation of our method in practice.
\section{Concluding discussion}\label{sec6}

In this paper, we considered maximum likelihood estimation
for the event time distribution function, $F$, at a grid point in the
current status model with i.i.d. data
and observation times lying on a regular grid. The spacing of the grid
$\delta$ was specified as $cn^{-\gamma}$
for constants $c>0$ and $0 < \gamma\leq1$ in order to incorporate
situations where
there are systematic ties in observation times, and the number of
distinct observation times can increase
with the sample size. The asymptotic properties of the NPMLE were shown
to depend on the order of the grid resolution $\gamma$
and an adaptive procedure, which circumvents the estimation of the
unknown~$\gamma$ and~$c$, was proposed for the construction of
asymptotically correct confidence intervals for the value of $F$ at a
grid-point of interest. We conclude with a description of alternative
methods for inference in this
problem and potential directions for future
research.\vspace*{8pt}

\textit{Likelihood ratio based inference}: An alternative to the Wald-type
adaptive confidence intervals proposed in this paper would be to use
those obtained via \mbox{likelihood} ratio inversion. More specifically, one
could consider testing the null hypothesis $H_0$ that $F(t_l) = \theta
_l$ versus its complement using the likelihood ratio statistics (LRS).
When the null hypothesis is true, the LRS converges weakly to $\chi
_1^2$ in the limit for $\gamma< 1/3$, to $\mathbb{D}$, the\vadjust{\goodbreak}
parameter-free limit discovered by \citet{Banerjee2001} for $\gamma>
1/3$ and a discrete analog of $\mathbb{D}$ depending on $c,\alpha
,\beta$, say $\mathcal{M}_{c,\alpha,\beta}$, that can be written in
terms of slopes of unconstrained and appropriately constrained convex
minorants of the process $\mathcal{P}_{c,\alpha,\beta}$ for $\gamma
= 1/3$. Thus, one obtains a boundary distribution for the likelihood
ratio statistic as well, and a phenomenon similar to that observed in
Section \ref{sec4} transpires, with the boundary distribution converging to
$\chi_1^2$ as $c \rightarrow\infty$ and to that of $\mathbb{D}$ as
$c \rightarrow0$. An adaptive procedure, which performs an inversion
by calibrating the likelihood ratio statistics for testing a family of
null hypotheses of the form $F(t_l) = \theta$ for varying $\theta$,
using the quantiles of $\mathcal{M}_{\hat{c},\hat{\alpha},\hat
{\beta}}$, can also be developed but is computationally more
burdensome than the Wald-type intervals. See \citet{tangbankos10} for
the details.\vspace*{8pt}

\textit{Smoothed estimators}: We recall that all our results have been
developed under minimal smoothness assumptions on $F$: throughout the
paper, we assume $F$ to be once continuously differentiable with a
nonvanishing derivative around $x_0$. We used the NPMLE to make
inference on $F$ since it can be computed without specifying
bandwidths; furthermore, under our minimal assumptions, its pointwise
rate of convergence when $\gamma> 1/3$ or when the observation times
arise from a continuous distribution cannot be bettered by a smoothed
estimator. However, if one makes the assumption of a second derivative
at $x_0$, the kernel-smoothed NPMLE (and related variants) can achieve
a convergence rate of $n^{2/5}$ (which is faster than the rate of the
NPMLE) using a bandwidth of order $n^{-1/5}$. See \citet
{Groeneboom2010} where these results are developed and also an earlier
paper due to \citet{Mammen1991} dealing with monotone regression. In
such a~situation, one could envisage using a~smoothed version of the
NPMLE in this problem with a bandwidth larger than the resolution of
the grid, and it is conceivable that an adaptive procedure could be
developed along these lines. While this is certainly an interesting and
important topic for further exploration, it is outside the scope of
this work, not least owing to the fact that the assumptions underlying
such a procedure are different (two derivatives as opposed to one) than
those in this paper.\vspace*{8pt}

\textit{Further possibilities}: The results in this paper reveal some
new directions for future research. As touched upon in the
\hyperref[intro]{Introduction}, some recent related work by
\citet{Maathuis2010} deals with the estimation of
\textit{competing risks current status data under finitely many risks}
with finitely many discrete (or grouped) observation times. A natural
question of interest, then, is what happens if the observation times in
their paper are supported on grids of increasing size as considered in
this paper for simple current status data. We suspect that a similar
adaptive procedure relying on a boundary phenomenon at $\gamma= 1/3$
can also be developed in this case. Furthermore, one could consider the
problem of grouped current status data\vadjust{\goodbreak} (with and without the element of
competing risks), where the observation times are not exactly known but
grouped into bins. Based on communications with us and preliminary
versions of this paper, \citet {Maathuis2010} conjecture that for
grouped current status data \textit{without competing risks}, one may
expect findings similar to those in this paper, depending on whether
the number of groups increases at rate $n^{1/3}$ or at a faster/slower
rate and it would not be unreasonable to expect a similar thing to
happen for grouped current status data with finitely many competing
risks. In fact, an adaptive inference procedure very similar to that in
this paper should also work for the problem treated in
\citet{Zhang2001} and allow inference for the decreasing density
of interest without needing to know the rate of growth of the bins.

It is also fairly clear that the adaptive inference scheme proposed in
this paper will apply to monotone regression models with discrete
covariates in general. In particular, the very general conditionally
parametric response models studied in \citet{Banerjee2007d} under the
assumption of a continuous covariate can be handled for the discrete
covariate case as well by adapting the methods of this paper.
Furthermore, similar adaptive inference in more complex forms of
interval censoring, like Case-2 censoring or mixed-case censoring [see,
e.g., \citet{Sen2006} and
\citet{Schick2000}], should also be possible in situations where the
multiple observation times are discrete-valued. Finally, we conjecture
that phenomena similar to those revealed in this paper will appear in
nonparametric regression problems with grid-supported covariates under
more complex shape constraints (like convexity, e.g.), though the
boundary value of $\gamma$ as well as the nature of the nonstandard
limits will be different and will depend on the ``order'' of the shape
constraint. This will also be a topic of future research.

%In the
%simple signal plus nosie convex regression model, heuristic
%calculations reveal indicate that grid resolutions with $\gamma< 1/5$
%will
%produce normal limits for least squares estimates and transition to a
%nonstandard limit will probably happen at $\gamma= 1/5$.
%though the boundary value for $\gamma$ would change
%and so would the asymptotic distributions for large $\gamma$.

%All the derivations and results in this paper with no modification can
%be applied to
%the general binary regression model with an increasing probability
%function,
%whose range is a subset of $[0,1]$.
%The assumption that the grid points are equally spaced
%can be weakened slightly but
%each distance between two successive grid points
%should lie between $c_{L}n^{-\gamma}$ and $c_{U}n^{-\gamma}$,
%where $c_{L}<c_{U}$ are two positive real constants.
%Similar asymptotic results
%of the three cases on $\gamma$ can be established
%for other isotonic regression problems under different models,
%which involve design densities.

\begin{appendix}\label{app}
%s7 #&#
\section*{Appendix: Proofs}

\begin{pf*}{Proof of Theorem
\ref{thmrelationshipfrom-one-third-to-lessWald}}
For $k\in\mathbb{Z}$, let
\[
\tilde h(k) = \alpha\sqrt{c} W(ck)+\beta c^{5/2} k(1+k),\qquad
h(k) = \alpha c W(k)+\beta c^{5/2} k(1+k).
\]
Then, we have
$
\{\tilde h(k), k\in\mathbb{Z}\}
\stackrel{d}{=}
\{h(k), k\in\mathbb{Z}\}
$.
Thus,
\[
\sqrt{c}
\mathcal{S}_{c}
\stackrel{d}{=}
\operatorname{LS}\circ \operatorname{GCM}
\{( ck, h(k) ), k\in\mathbb{Z} \}(0).
\]
Define $\tilde{\mathcal S}_c = \sqrt{c} \mathcal{S}_c$.
Denote
\begin{eqnarray*}
A_{c} &=&
\biggl\{\frac{h(k)}{ck} < \frac{h(k+1)}{c(k+1)}, k = 1,2,\ldots\biggr\}, \\
B_{c} &=&
\biggl\{\frac{h(-(k-1))}{c(k-1)} < \frac{h(-k)}{ck}, k = 2, 3, \ldots
\biggr\},\\
C_{c} &=&
\biggl\{\frac{h(1)}{c} > \frac{-h(-1)}{c} \biggr\}.
\end{eqnarray*}
Then, for $\omega\in A_{c}B_{c}C_{c}$,
it is easy to see $\tilde{\mathcal{S}}_{c}=-\alpha W(-1)$.
We will show in Lem\-ma~\ref{Lemmarelationshipcasetwotocaseone},
$A_{c}B_{c}C_{c}\stackrel{P}{\rightarrow}1$.
Thus,
$
\tilde{\mathcal{S}}_{c}
= \tilde{\mathcal{S}}_{c}A_{c}B_{c}C_{c}
+ \tilde{\mathcal{S}}_{c}(1-A_{c}B_{c}C_{c})
\stackrel{d}{\rightarrow}
-\alpha W(-1)
\stackrel{d}{=}
\alpha Z
$,
with $Z\sim N(0,1)$.
Therefore, $\sqrt{c}\mathcal{S}_{c}\stackrel{d}{\rightarrow}\alpha Z$.
\end{pf*}
%
%le7.1 #&#
\begin{lem}\label{Lemmarelationshipcasetwotocaseone}
Each of $A_{c}$, $B_{c}$ and $C_{c}$ in the proof
of Theorem \ref{thmrelationshipfrom-one-third-to-lessWald}
converges to 1 in probability.
\end{lem}
\begin{pf}
It is easy to show $C_{c}$ converges
to 1 in probability. The argument that $A_{c}$ converges to one in
probability is similar to that for $B_c$, and
only the former is established here.
In order to show $P(A_{c})\rightarrow1$,
it suffices to show $P(A_{c}^{c})\rightarrow0$.
We have, for each $k\in\mathbb{Z}$,
\begin{eqnarray*}
& & P\biggl(\frac{h(k)}{ck} \geq\frac{h(k+1)}{c(k+1)} \biggr)\\
%&=&P (\frac{\alpha c W(k)}{k} + \beta c^{5/2}(k+1) \geq\frac{
&&\qquad=P \biggl(\frac{\alpha W(k)}{k} + \beta c^{3/2}(k+1) \geq\frac
{\alpha W(k+1)}{k+1} + \beta c^{3/2}(k+2) \biggr)\\
&&\qquad = P\biggl(\alpha\biggl[\frac{W(k)}{k} - \frac{W(k+1)}{k+1} \biggr]
\geq\beta c^{3/2} \biggr)\\
%&=& P(\alpha N(0, k(k+1)) \geq\beta c^{3/2} k(k+1) )\\
&&\qquad = P\bigl( N(0,1) \geq\alpha^{-1}\beta c^{3/2}\sqrt{k(k+1)} \bigr)\\
&&\qquad\leq 2^{-1}\exp\{- 2^{-1}\alpha^{-2}\beta^{2}c^{3}k(k+1) \}
\end{eqnarray*}
using the fact that $W(k)/k - W(k+1)/(k+1) \sim N(0,(k (k+1))^{-1})$
and the inequality $P(N(0,1)>x)\leq2^{-1}\exp\{(-2^{-1}x^{2})\}$ for
$x\geq0$
[see, e.g., $\langle2\rangle$ on pa\-ge~317 of \citet{Pollard2002}].
Then, we have
\begin{eqnarray*}
P(A_{c}^{c})
&\leq&\sum_{k=1}^{\infty}P\biggl(\frac{h(k)}{ck} \geq\frac
{h(k+1)}{c(k+1)} \biggr)
\leq\sum_{k=1}^{\infty}2^{-1}\exp\{- 2^{-1}\alpha^{-2}\beta
^{2}c^{3}k^{2} \}\\
&\leq&2^{-1}\int_{0}^{\infty}
\exp\{-2^{-1}\alpha^{-2}\beta^{2}c^{3}x^{2} \} \,dx
= \bigl(\sqrt{2\pi}/4\bigr)\alpha\beta^{-1}c^{-3/2}
\rightarrow0
\end{eqnarray*}
as $c\rightarrow\infty$.
Thus, $P(A_{c})\rightarrow1$,
which completes the proof.
\end{pf}
\begin{pf*}{Proof of Theorem \ref
{thmrelationshipfrom-one-third-to-moreWald}}
We want to show that
$ \mathcal{S}_{c} \stackrel{d}{\rightarrow}
g_{\alpha,\beta}(0),
$ as \mbox{$c\rightarrow0$},
where
$
g_{\alpha,\beta}(0)
= \operatorname{LS}\circ \operatorname{GCM} \{ X_{\alpha,\beta}\}(0)
= \operatorname{LS}\circ \operatorname{GCM} \{ X_{\alpha,\beta}(t)\dvtx t\in\mathbb{R}\}(0)
$ and\break
\mbox{$
\mathcal{S}_{c}
= \operatorname{LS}\circ \operatorname{GCM}\{\mathcal{P}_{c}\}(0)
= \operatorname{LS}\circ \operatorname{GCM}\{\mathcal{P}_{c}(k)\dvtx k\in\mathbb{Z}\}(0)
$}.
%Since
% \mathcal{S}_{c}
% & = \operatorname{LS}\circ GCM\{\mathcal{P}_{c}(k): k\in\mathbb{Z}\}(0) \\
% & = \operatorname{LS}\circ GCM\{(ck, \alpha W(ck) + \beta c^{2}k(1+k)): k\in
% & = \operatorname{LS}\circ GCM\{(ck, \alpha W(ck) + \beta(ck)^{2} + \beta c (ck)): k
% & = \operatorname{LS}\circ GCM\{(ck, \alpha W(ck) + \beta(ck)^{2}): k\in\mathbb{Z}
% & + \operatorname{LS}\circ GCM\{(ck, \alpha W(ck) + \beta c (ck)): k\in\mathbb{Z}
% & = \operatorname{LS}\circ GCM\{(ck, \alpha W(ck) + \beta(ck)^{2}): k\in\mathbb{Z}
% & =: \operatorname{LS}\circ GCM\{\mathcal{P}_{c}': k\in\mathbb{Z}\}(0) + \beta c \\
% & =: \mathcal{S}_{c}' + \beta c,
Since $\mathcal{S}_{c}=\mathcal{S}_{c}' + \beta c$,\break
where $\mathcal{S}_{c}'=\operatorname{LS}\circ \operatorname{GCM}\{\mathcal{P}_{c}'\dvtx k\in\mathbb
{Z}\}(0)$
and
$\mathcal{P}_{c}'=\{(ck, \alpha W(ck) + \beta(ck)^{2})\dvtx\break k\in\mathbb
{Z}\}$,\vspace*{-1pt}
it is sufficient to show
$
\mathcal{S}_{c}' \stackrel{d}{\rightarrow}
g_{\alpha,\beta}(0)
$
as $c\rightarrow0$.
To make the notation simple and without causing confusion,
in the following
we still use~$\mathcal{P}_{c}$ and $\mathcal{S}_{c}$
to denote
$\mathcal{P}_{c}'$ and $\mathcal{S}_{c}'$. Also, it will be useful to
think of $\mathcal{P}_c$ as
a continuous process on $\RR$ formed by linearly interpolating the
points $\{ck, \mathcal{P}_{2,c}(ck)\dvtx k \in\mathbb{Z}\}$, where
$\mathcal{P}_{2,c}(ck) = \alpha W(ck) + \beta (ck)^2 = X_{\alpha
,\beta}(ck)$. Note that\vadjust{\goodbreak} viewing $\mathcal{P}_c$ in this way keeps the GCM
unaltered, that is, the GCM of this continuous linear interpolated
version is the same as that of the set of points $\{ck, \mathcal
{P}_{2,c}(ck)\dvtx k \in\mathbb{Z}\}$,
and the slope-changing points of this piece-wise linear GCM are still
grid-points of the form $ck$.

Let $L$ and $U$
be the largest negative and smallest nonnegative
$x$-axis coordinates of the slope changing points of the GCM of
$X_{\alpha,\beta}$. Similarly,
let~$L_{c}$ and $U_{c}$
be the largest negative and smallest nonnegative
$x$-axis coordinates
of the slope changing points of the GCM of $\mathcal{P}_{c}$. For $K >
0$, define $g_{\alpha,\beta}^K(0) = \operatorname{LS} \circ \operatorname{GCM}\{X_{\alpha,\beta
}(t)\dvtx t \in[-K,K]\}(0)$ and $\mathcal{S}_c^{K} = \operatorname{LS} \circ \operatorname{GCM}\{
\mathcal{P}_c(t) \dvtx\break t \in[-K,K]\}(0)$.

We will show that,
given $\varepsilon> 0$,
there exist $M_{\varepsilon} > 0$ and $c(\varepsilon)$
such that (a) for all $0 < c < c(\varepsilon)$, $P(\mathcal
{S}_c^{M_{\varepsilon}}
\ne\mathcal{S}_c) < \varepsilon$ and (b) $P(g_{\alpha,\beta
}^{M_{\varepsilon}}(0) \ne g_{\alpha,\beta}(0)) < \varepsilon$. These
immediately imply
that both Fact \ref{fact1}: $\lim_{\varepsilon\rightarrow0} \limsup_{c
\rightarrow0} P(\mathcal{S}_{c}^{M_{\varepsilon}} \not= \mathcal
{S}_{c}) = 0$ and
Fact~\ref{fact2}: $\lim_{\varepsilon\rightarrow0}P(g_{\alpha,\beta
}^{M_{\varepsilon}}(0) \not= g_{\alpha,\beta}(0)) = 0$ hold. We then
show that
Fact \ref{fact3}: for\vspace*{-1pt} each $\varepsilon>0$, $\mathcal{S}_{c}^{M_{\varepsilon}}
\stackrel{d}{\rightarrow} g_{\alpha,\beta}^{M_{\varepsilon}}(0)$
holds as well. Then, by Lemma \ref{LemmaTruncation},
we have the conclusion $\mathcal{S}_{c}\stackrel{d}{\rightarrow
}g_{\alpha,\beta}(0)$.
Figure \ref{fgdrawingLimit-Localization}
illustrates the following argument.

%f3 #&#
\begin{figure}

\includegraphics{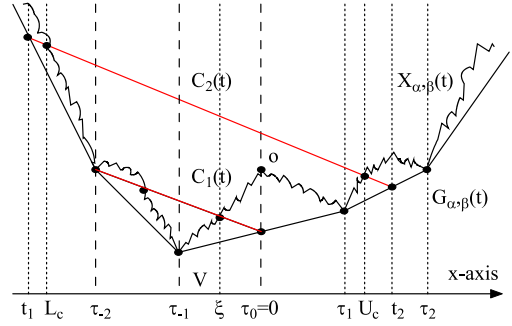}

\caption{An illustration for showing $\{L_{c}\}$ is $O_{P}(1)$
in the proof of
Theorem \protect\ref{thmrelationshipfrom-one-third-to-moreWald}.}
\label{fgdrawingLimit-Localization}
\end{figure}

Let $\tau_{-2} < \tau_{-1} < \tau_1 < \tau_2$ be four consecutive
slope changing points of $G_{\alpha,\beta} = \operatorname{GCM}\{X_{\alpha,\beta}\}$
with $\tau_{-1}$ denoting the first slope changing point to the left
of 0 and $\tau_1$ the first slope changing point to the right.
%Similarly define $\tau_{-2}$ and $\tau_{2}$.
Since $\tau_{-2}$
and $\tau_2$ are $O_{P}(1)$, given $\varepsilon> 0$,
there exists $M_{\varepsilon} > 0$ such that $P(-M_{\varepsilon} < \tau
_{-2} < \tau_{2} < M_{\varepsilon})
> 1 - \varepsilon/4$. Note that the event
$\{g_{\alpha,\beta}^{M_{\varepsilon}}(0) = g_{\alpha,\beta}(0)\}
\subset\{-M_{\varepsilon} < \tau_{-2} < \tau_{2} < M_{\varepsilon}\}$,
and it follows that
$P(g_{\alpha,\beta}^{M_{\varepsilon}}(0) \not= g_{\alpha,\beta}(0))
< \varepsilon/4
< \varepsilon$.
Thus, (b) holds.\vspace*{1pt}

Next, consider the chord $C_1(t)$ joining $(0, G_{\alpha,\beta}(0))$
and $(\tau_{-2}, G_{\alpha,\beta}(\tau_{-2}))$.
By the convexity of $G_{\alpha,\beta}$ over $[\tau_{-2},0]$
and $\tau_{-1}\in(\tau_{-2},0)$ being a slope changing point,
$X_{\alpha,\beta}(\tau_{-1}) = G_{\alpha,\beta}(\tau_{-1}) <
C(\tau_{-1})$.
But $C_1(0) = G_{\alpha,\beta}(0) < X_{\alpha,\beta}(0)$, and
it follows by the intermediate value theorem that
$\xi= \inf_{\tau_{-1} < t < 0} \{t\dvtx\break X_{\alpha,\beta}(t) = C_1(t)\}
$ is well defined
(since the set in question is nonempty),
$ \tau_{-1} < \xi< 0$,
$C_1(\xi) = X_{\alpha,\beta}(\xi)$ and
on $[\tau_{-1}, \xi)$, $X_{\alpha,\beta}(t) < C_1(t)$.
Let $V = \xi- \tau_{-1}$.
Since $V$ is a continuous and
positive random variable,
there exists $\delta(\varepsilon) > 0$
such that $P(V > \delta(\varepsilon)) \geq1 - \varepsilon/4$.
Then, the event $E_{\varepsilon} = \{V > \delta(\varepsilon)\} \cap\{
-M_{\varepsilon} < \tau_{-2}\}$
has probability larger than $1 - \varepsilon/2$.
For any $c < c(\varepsilon) =: \delta(\varepsilon)$,
we claim that $L_c \geq\tau_{-2}$ on the event $E_{\varepsilon}$,
and the argument for this follows below.

If $L_c < \tau_{-2}$, consider the chord $C_2(t)$
connecting two points $(L_{c}, \mathcal{P}_{2,c}(L_{c}))$ and $(U_c,
\mathcal{P}_{2,c}(U_{c}))$.
This chord must lie strictly above the chord $\{C_1(t)\dvtx \tau_{-1} \leq
t \leq0\}$
since it can be viewed as a restriction of a chord connecting two points
$(t_1, G_{\alpha,\beta}(t_1))$ and $(t_2, G_{\alpha,\beta}(t_2))$
with $t_1\leq L_c < \tau_{-1} < 0 \leq U_c \leq t_2$. It then follows
that all points of the form
$\{ck, \mathcal{P}_{2,c}(ck) = X_{\alpha,\beta}(ck)\dvtx ck \in[L_c,
U_c]\}$
must lie above $C_2(t)$. But there is at least one $ck^{\star}$ with
$\tau_{-1} < ck^{\star} < \xi$ and
such that $X_{\alpha,\beta}(ck^{\star}) < C_1(ck^{\star}) <
C_2(ck^{\star})$,
which furnishes a contradiction.

We conclude that for any $c < c(\varepsilon)$,
$P(-M_{\varepsilon} < L_c) > 1 - \varepsilon/2$.
A similar argument to the right-hand side of 0 shows that for the same $c$'s
(by the symmetry of two-sided Brownian motion about the origin), $P(U_c
< M_{\varepsilon}) > 1 - \varepsilon/2$. Hence $P(-M_{\varepsilon} < L_c <
U_c <
M_{\varepsilon}) > 1 - \varepsilon$. On this event, clearly $\mathcal
{S}_c^{M_{\varepsilon}} = \mathcal{S}_c$, and it follows that for all $c
< c(\varepsilon)$,
$P(\mathcal{S}_c^{M_{\varepsilon}} \not= \mathcal{S}_c) < \varepsilon$.
Thus, (a) also holds and Facts \ref{fact1} and~\ref{fact2} are established.

It remains to establish Fact \ref{fact3}. This follows easily.
%Consider the set $[\pm M_{\varepsilon}]$.
For almost every $\omega$, $X_{\alpha,\beta}(t)$
is uniformly continuous on $[\pm2 M_{\varepsilon}]$.
It follows by elementary analysis that (for almost every $\omega$)
on $[\pm M_{\varepsilon}]$, the process $\mathcal{P}_c$,
being the linear interpolant of the points
$\{ck, X_{\alpha,\beta}(ck)\dvtx -M_{\varepsilon} \leq ck \leq M_{\varepsilon
}\}
\cup\{ (-M_\varepsilon, \mathcal{P}_{2c}(-M_\varepsilon)),\break
(M_\varepsilon, \mathcal{P}_{2,c}(M_\varepsilon))\}$,
converges uniformly to
$X_{\alpha,\beta}$ as $c \rightarrow0$. Thus, the left slope of the
GCM of $\{\mathcal{P}_c(t)\dvtx t \in[\pm M_{\varepsilon}]\}$, which is precisely
$\mathcal{S}_c^{M_{\varepsilon}}$, converges to $g_{\alpha,\beta
}^{M_{\varepsilon}}(0)$
since the GCM of the restriction of $X_{\alpha,\beta}$ to
$[\pm M_{\varepsilon}]$ is almost surely differentiable at 0; see, for
example, the Lemma on page 330 of \citet{Robertson1988}
for a justification of this convergence.
\end{pf*}
\end{appendix}

%
%%\newpage
%%\input{Appendix_Case_Two_033_1}
%
%%\newpage
%
%%\newpage

\section*{Acknowledgments}

We would like to thank Professors Jack Kalbfleisch and Nick Jewell for
bringing this problem to our attention. The first author would also
like to thank Professor George Michailidis for partial financial
support while he was involved with the project.

\begin{supplement}[id=suppA]
\stitle{More proofs for the current paper ``Likelihood based
inference for current status data on a grid: A boundary phenomenon and
an adaptive inference procedure''}
\slink[doi]{10.1214/11-AOS942SUPP} %[doi,text={...}] - jei reikia
%suskaldyti doi
\sdatatype{.pdf}
\sfilename{aos942\_supp.pdf}
\sdescription{The supplementary material contains the details of the
proofs of several theorems and lemmas in
Sections \ref{secasymptoticsgammaless033} and \ref
{secasymptoticsgammaequal033} of this paper.}
\end{supplement}

% imsref loaded by lrinkeviciute, 2012-01-20 13:25:03
% imsref loaded by lrinkeviciute, 2012-01-20 13:31:22

\printaddresses


\begin{thebibliography}{30}
% BibTex style file: ims.bst, 2011-05-30
% Default style options (sort=0,type=number).
% Used options (sort=1,type=nameyear).

%b1 ###
\bibitem[\protect\citeauthoryear{Abrevaya and Huang}{2005}]{Abrevaya2005}
\begin{barticle}[mr]
\bauthor{\bsnm{Abrevaya},~\bfnm{Jason}\binits{J.}} \AND
  \bauthor{\bsnm{Huang},~\bfnm{Jian}\binits{J.}}
(\byear{2005}).
\btitle{On the bootstrap of the maximum score estimator}.
\bjournal{Econometrica}
\bvolume{73}
\bpages{1175--1204}.
\bid{doi={10.1111/j.1468-0262.2005.00613.x}, issn={0012-9682}, mr={2149245}}
\bptok{imsref}%
\end{barticle}
\endbibitem

%b2 ###
\bibitem[\protect\citeauthoryear{Anevski and H{\"o}ssjer}{2006}]{Anevski2006}
\begin{barticle}[mr]
\bauthor{\bsnm{Anevski},~\bfnm{D.}\binits{D.}} \AND
  \bauthor{\bsnm{H{\"o}ssjer},~\bfnm{O.}\binits{O.}}
(\byear{2006}).
\btitle{A general asymptotic scheme for inference under order restrictions}.
\bjournal{Ann. Statist.}
\bvolume{34}
\bpages{1874--1930}.
\bid{doi={10.1214/009053606000000443}, issn={0090-5364}, mr={2283721}}
\bptok{imsref}%
\end{barticle}
\endbibitem

%b3 ###
\bibitem[\protect\citeauthoryear{Banerjee}{2007}]{Banerjee2007d}
\begin{barticle}[mr]
\bauthor{\bsnm{Banerjee},~\bfnm{Moulinath}\binits{M.}}
(\byear{2007}).
\btitle{Likelihood based inference for monotone response models}.
\bjournal{Ann. Statist.}
\bvolume{35}
\bpages{931--956}.
\bid{doi={10.1214/009053606000001578}, issn={0090-5364}, mr={2341693}}
\bptok{imsref}%
\end{barticle}
\endbibitem

%b4 ###
\bibitem[\protect\citeauthoryear{Banerjee and Wellner}{2001}]{Banerjee2001}
\begin{barticle}[mr]
\bauthor{\bsnm{Banerjee},~\bfnm{Moulinath}\binits{M.}} \AND
  \bauthor{\bsnm{Wellner},~\bfnm{Jon~A.}\binits{J.~A.}}
(\byear{2001}).
\btitle{Likelihood ratio tests for monotone functions}.
\bjournal{Ann. Statist.}
\bvolume{29}
\bpages{1699--1731}.
\bid{doi={10.1214/aos/1015345959}, issn={0090-5364}, mr={1891743}}
\bptok{imsref}%
\end{barticle}
\endbibitem

%b5 ###
\bibitem[\protect\citeauthoryear{Banerjee and Wellner}{2005}]{Banerjee2005}
\begin{barticle}[mr]
\bauthor{\bsnm{Banerjee},~\bfnm{Moulinath}\binits{M.}} \AND
  \bauthor{\bsnm{Wellner},~\bfnm{Jon~A.}\binits{J.~A.}}
(\byear{2005}).
\btitle{Confidence intervals for current status data}.
\bjournal{Scand. J.~Stat.}
\bvolume{32}
\bpages{405--424}.
\bid{doi={10.1111/j.1467-9469.2005.00454.x}, issn={0303-6898}, mr={2204627}}
\bptok{imsref}%
\end{barticle}
\endbibitem

%b6 ###
\bibitem[\protect\citeauthoryear{Brunk}{1970}]{Brunk1970}
\begin{bincollection}[mr]
\bauthor{\bsnm{Brunk},~\bfnm{H.~D.}\binits{H.~D.}}
(\byear{1970}).
\btitle{Estimation of isotonic regression}.
In \bbooktitle{Nonparametric {T}echniques in {S}tatistical {I}nference ({P}roc.
  {S}ympos., {I}ndiana {U}niv., {B}loomington, {I}nd., 1969)}
\bpages{177--197}.
\bpublisher{Cambridge Univ. Press}, \baddress{London}.
\bid{mr={0277070}}
\bptok{imsref}%
\end{bincollection}
\endbibitem

%b7 ###
\bibitem[\protect\citeauthoryear{Groeneboom, Jongbloed and
  Witte}{2010}]{Groeneboom2010}
\begin{barticle}[mr]
\bauthor{\bsnm{Groeneboom},~\bfnm{Piet}\binits{P.}},
  \bauthor{\bsnm{Jongbloed},~\bfnm{Geurt}\binits{G.}} \AND
  \bauthor{\bsnm{Witte},~\bfnm{Birgit~I.}\binits{B.~I.}}
(\byear{2010}).
\btitle{Maximum smoothed likelihood estimation and smoothed maximum likelihood
  estimation in the current status model}.
\bjournal{Ann. Statist.}
\bvolume{38}
\bpages{352--387}.
\bid{doi={10.1214/09-AOS721}, issn={0090-5364}, mr={2589325}}
\bptok{imsref}%
\end{barticle}
\endbibitem

%b8 ###
\bibitem[\protect\citeauthoryear{Groeneboom and Wellner}{1992}]{Groeneboom1992}
\begin{bbook}[mr]
\bauthor{\bsnm{Groeneboom},~\bfnm{Piet}\binits{P.}} \AND
  \bauthor{\bsnm{Wellner},~\bfnm{Jon~A.}\binits{J.~A.}}
(\byear{1992}).
\btitle{Information Bounds and Nonparametric Maximum Likelihood Estimation}.
\bseries{DMV Seminar}
\bvolume{19}.
\bpublisher{Birkh\"auser}, \baddress{Basel}.
\bid{mr={1180321}}
\bptok{imsref}%
\end{bbook}
\endbibitem

%b9 ###
\bibitem[\protect\citeauthoryear{Keiding et~al.}{1996}]{Keiding1996}
\begin{barticle}[author]
\bauthor{\bsnm{Keiding},~\bfnm{Neils}\binits{N.}},
  \bauthor{\bsnm{Begtrup},~\bfnm{Kamilla}\binits{K.}},
  \bauthor{\bsnm{Scheike},~\bfnm{Thomas~H.}\binits{T.~H.}} \AND
  \bauthor{\bsnm{Hasibeder},~\bfnm{Gunther}\binits{G.}}
(\byear{1996}).
\btitle{Estimation from current status data in continuous time}.
\bjournal{Lifetime Data Anal.}
\bvolume{2}
\bpages{119--129}.
\bptok{imsref}%
\end{barticle}
\endbibitem

%b10 ###
\bibitem[\protect\citeauthoryear{Kiefer and Wolfowitz}{1976}]{Kiefer1976}
\begin{barticle}[mr]
\bauthor{\bsnm{Kiefer},~\bfnm{J.}\binits{J.}} \AND
  \bauthor{\bsnm{Wolfowitz},~\bfnm{J.}\binits{J.}}
(\byear{1976}).
\btitle{Asymptotically minimax estimation of concave and convex distribution
  functions}.
\bjournal{Z. Wahrsch. Verw. Gebiete}
\bvolume{34}
\bpages{73--85}.
\bid{mr={0397974}}
\bptok{imsref}%
\end{barticle}
\endbibitem

%b11 ###
\bibitem[\protect\citeauthoryear{Kim and Pollard}{1990}]{Kim1990}
\begin{barticle}[mr]
\bauthor{\bsnm{Kim},~\bfnm{JeanKyung}\binits{J.}} \AND
  \bauthor{\bsnm{Pollard},~\bfnm{David}\binits{D.}}
(\byear{1990}).
\btitle{Cube root asymptotics}.
\bjournal{Ann. Statist.}
\bvolume{18}
\bpages{191--219}.
\bid{doi={10.1214/aos/1176347498}, issn={0090-5364}, mr={1041391}}
\bptok{imsref}%
\end{barticle}
\endbibitem

%b12 ###
\bibitem[\protect\citeauthoryear{Kosorok}{2008}]{Kosorok2008a}
\begin{bincollection}[author]
\bauthor{\bsnm{Kosorok},~\bfnm{Michael~R.}\binits{M.~R.}}
(\byear{2008}).
\btitle{Bootstrapping the Grenander estimator}.
In \bbooktitle{Beyond Parametrics in Interdisciplinary Research: Festschrift in
  Honor of Professor Pranab K. Sen}
\bpages{282--292}.
\bpublisher{IMS}, \baddress{Hayward, CA}.
\bptok{imsref}%
\end{bincollection}
\endbibitem

%b13 ###
\bibitem[\protect\citeauthoryear{Leurgans}{1982}]{Leurgans1982}
\begin{barticle}[mr]
\bauthor{\bsnm{Leurgans},~\bfnm{Sue}\binits{S.}}
(\byear{1982}).
\btitle{Asymptotic distributions of slope-of-greatest-convex-minorant
  estimators}.
\bjournal{Ann. Statist.}
\bvolume{10}
\bpages{287--296}.
\bid{issn={0090-5364}, mr={0642740}}
\bptok{imsref}%
\end{barticle}
\endbibitem

%b14 ###
\bibitem[\protect\citeauthoryear{Maathuis and Hudgens}{2011}]{Maathuis2010}
\begin{barticle}[mr]
\bauthor{\bsnm{Maathuis},~\bfnm{M.~H.}\binits{M.~H.}} \AND
  \bauthor{\bsnm{Hudgens},~\bfnm{M.~G.}\binits{M.~G.}}
(\byear{2011}).
\btitle{Nonparametric inference for competing risks current status data with
  continuous, discrete or grouped observation times}.
\bjournal{Biometrika}
\bvolume{98}
\bpages{325--340}.
\bid{doi={10.1093/biomet/asq083}, issn={0006-3444}, mr={2806431}}
\bptok{imsref}%
\end{barticle}
\endbibitem

%b15 ###
\bibitem[\protect\citeauthoryear{Mammen}{1991}]{Mammen1991}
\begin{barticle}[mr]
\bauthor{\bsnm{Mammen},~\bfnm{Enno}\binits{E.}}
(\byear{1991}).
\btitle{Estimating a smooth monotone regression function}.
\bjournal{Ann. Statist.}
\bvolume{19}
\bpages{724--740}.
\bid{doi={10.1214/aos/1176348117}, issn={0090-5364}, mr={1105841}}
\bptok{imsref}%
\end{barticle}
\endbibitem

%b16 ###
\bibitem[\protect\citeauthoryear{Pollard}{2002}]{Pollard2002}
\begin{bbook}[author]
\bauthor{\bsnm{Pollard},~\bfnm{David}\binits{D.}}
(\byear{2002}).
\btitle{A User's Guide to Measure Theoretic Probability}.
\bpublisher{Cambridge Univ. Press}, \baddress{Cambridge}.
\bptok{imsref}%
\end{bbook}
\endbibitem

%b17 ###
\bibitem[\protect\citeauthoryear{Prakasa~Rao}{1969}]{Rao1969}
\begin{barticle}[mr]
\bauthor{\bsnm{Prakasa~Rao},~\bfnm{B.~L.~S.}\binits{B.~L.~S.}}
(\byear{1969}).
\btitle{Estkmation of a unimodal density}.
\bjournal{Sankhy\=a Ser. A}
\bvolume{31}
\bpages{23--36}.
\bid{issn={0581-572X}, mr={0267677}}
\bptok{imsref}%
\end{barticle}
\endbibitem

%b18 ###
\bibitem[\protect\citeauthoryear{Robertson, Wright and
  Dykstra}{1988}]{Robertson1988}
\begin{bbook}[mr]
\bauthor{\bsnm{Robertson},~\bfnm{Tim}\binits{T.}},
  \bauthor{\bsnm{Wright},~\bfnm{F.~T.}\binits{F.~T.}} \AND
  \bauthor{\bsnm{Dykstra},~\bfnm{R.~L.}\binits{R.~L.}}
(\byear{1988}).
\btitle{Order Restricted Statistical Inference}.
\bpublisher{Wiley}, \baddress{Chichester}.
\bid{mr={0961262}}
\bptok{imsref}%
\end{bbook}
\endbibitem

%b19 ###
\bibitem[\protect\citeauthoryear{Schick and Yu}{2000}]{Schick2000}
\begin{barticle}[mr]
\bauthor{\bsnm{Schick},~\bfnm{Anton}\binits{A.}} \AND
  \bauthor{\bsnm{Yu},~\bfnm{Qiqing}\binits{Q.}}
(\byear{2000}).
\btitle{Consistency of the {GMLE} with mixed case interval-censored data}.
\bjournal{Scand. J. Stat.}
\bvolume{27}
\bpages{45--55}.
\bid{doi={10.1111/1467-9469.00177}, issn={0303-6898}, mr={1774042}}
\bptok{imsref}%
\end{barticle}
\endbibitem

%b20 ###
\bibitem[\protect\citeauthoryear{Sen and Banerjee}{2007}]{Sen2006}
\begin{barticle}[mr]
\bauthor{\bsnm{Sen},~\bfnm{Bodhisattva}\binits{B.}} \AND
  \bauthor{\bsnm{Banerjee},~\bfnm{Moulinath}\binits{M.}}
(\byear{2007}).
\btitle{A pseudolikelihood method for analyzing interval censored data}.
\bjournal{Biometrika}
\bvolume{94}
\bpages{71--86}.
\bid{doi={10.1093/biomet/asm011}, issn={0006-3444}, mr={2307901}}
\bptok{imsref}%
\end{barticle}
\endbibitem

%b21 ###
\bibitem[\protect\citeauthoryear{Sen, Banerjee and Woodroofe}{2010}]{Sen2010}
\begin{barticle}[mr]
\bauthor{\bsnm{Sen},~\bfnm{Bodhisattva}\binits{B.}},
  \bauthor{\bsnm{Banerjee},~\bfnm{Moulinath}\binits{M.}} \AND
  \bauthor{\bsnm{Woodroofe},~\bfnm{Michael}\binits{M.}}
(\byear{2010}).
\btitle{Inconsistency of bootstrap: The {G}renander estimator}.
\bjournal{Ann. Statist.}
\bvolume{38}
\bpages{1953--1977}.
\bid{doi={10.1214/09-AOS777}, issn={0090-5364}, mr={2676880}}
\bptok{imsref}%
\end{barticle}
\endbibitem

%b22 ###
\bibitem[\protect\citeauthoryear{Tang, Banerjee and
  Kosorok}{2010}]{tangbankos10}
\begin{bmisc}[author]
\bauthor{\bsnm{Tang},~\bfnm{Runlong}\binits{R.}},
  \bauthor{\bsnm{Banerjee},~\bfnm{Moulinath}\binits{M.}} \AND
  \bauthor{\bsnm{Kosorok},~\bfnm{Michael~R.}\binits{M.~R.}}
(\byear{2010}).
\bhowpublished{Asymptotics for current status data under varying observation
  time sparsity. Available at
\href{http://www.stat.lsa.umich.edu/\textasciitilde moulib/csdgriddec23.pdf}%
    {www.stat.lsa.umich.edu/}
\href{http://www.stat.lsa.umich.edu/\textasciitilde moulib/csdgriddec23.pdf}%
    {\textasciitilde moulib/csdgriddec23.pdf}.}
\bptok{imsref}%
\end{bmisc}
\endbibitem

%b23 ###
\bibitem[\protect\citeauthoryear{Tang, Banerjee and Kosorok}{2011}]{Tang2011}
\begin{bmisc}[author]
\bauthor{\bsnm{Tang},~\bfnm{Runlong}\binits{R.}},
  \bauthor{\bsnm{Banerjee},~\bfnm{Moulinath}\binits{M.}} \AND
  \bauthor{\bsnm{Kosorok},~\bfnm{Michael~R.}\binits{M.~R.}}
(\byear{2011}).
\bhowpublished{Supplement to ``Likelihood based inference for current status
  data on a grid: A boundary phenomenon and an adaptive inference procedure.''
  DOI:\href{http://dx.doi.org/10.1214/11-AOS942SUPP}{10.1214/11-AOS942SUPP}.}
\bptok{imsref}%
\end{bmisc}
\endbibitem

%b24 ###
\bibitem[\protect\citeauthoryear{van~der Vaart}{1991}]{Vaart1991}
\begin{barticle}[mr]
\bauthor{\bparticle{van~der} \bsnm{Vaart},~\bfnm{Aad}\binits{A.}}
(\byear{1991}).
\btitle{On differentiable functionals}.
\bjournal{Ann. Statist.}
\bvolume{19}
\bpages{178--204}.
\bid{doi={10.1214/aos/1176347976}, issn={0090-5364}, mr={1091845}}
\bptok{imsref}%
\end{barticle}
\endbibitem

%b25 ###
\bibitem[\protect\citeauthoryear{van~der Vaart and van~der
  Laan}{2003}]{Vaart2003}
\begin{barticle}[mr]
\bauthor{\bparticle{van~der} \bsnm{Vaart},~\bfnm{Aad~W.}\binits{A.~W.}} \AND
  \bauthor{\bparticle{van~der} \bsnm{Laan},~\bfnm{Mark~J.}\binits{M.~J.}}
(\byear{2003}).
\btitle{Smooth estimation of a~monotone density}.
\bjournal{Statistics}
\bvolume{37}
\bpages{189--203}.
\bid{doi={10.1080/0233188031000124392}, issn={0233-1888}, mr={1986176}}
\bptok{imsref}%
\end{barticle}
\endbibitem

%b26 ###
\bibitem[\protect\citeauthoryear{Wellner and Zhang}{2000}]{Wellner2000}
\begin{barticle}[mr]
\bauthor{\bsnm{Wellner},~\bfnm{Jon~A.}\binits{J.~A.}} \AND
  \bauthor{\bsnm{Zhang},~\bfnm{Ying}\binits{Y.}}
(\byear{2000}).
\btitle{Two estimators of the mean of a counting process with panel count
  data}.
\bjournal{Ann. Statist.}
\bvolume{28}
\bpages{779--814}.
\bid{doi={10.1214/aos/1015951998}, issn={0090-5364}, mr={1792787}}
\bptok{imsref}%
\end{barticle}
\endbibitem

%b27 ###
\bibitem[\protect\citeauthoryear{Wright}{1981}]{Wright1981}
\begin{barticle}[mr]
\bauthor{\bsnm{Wright},~\bfnm{F.~T.}\binits{F.~T.}}
(\byear{1981}).
\btitle{The asymptotic behavior of monotone regression estimates}.
\bjournal{Ann. Statist.}
\bvolume{9}
\bpages{443--448}.
\bid{issn={0090-5364}, mr={0606630}}
\bptok{imsref}%
\end{barticle}
\endbibitem

%b28 ###
\bibitem[\protect\citeauthoryear{Yu et~al.}{1998}]{Yu1998}
\begin{barticle}[mr]
\bauthor{\bsnm{Yu},~\bfnm{Qiqing}\binits{Q.}},
  \bauthor{\bsnm{Schick},~\bfnm{Anton}\binits{A.}},
  \bauthor{\bsnm{Li},~\bfnm{Linxiong}\binits{L.}} \AND
  \bauthor{\bsnm{Wong},~\bfnm{George Y.~C.}\binits{G.~Y.~C.}}
(\byear{1998}).
\btitle{Asymptotic properties of the {GMLE} in the case {$1$}
  interval-censorship model with discrete inspection times}.
\bjournal{Canad. J. Statist.}
\bvolume{26}
\bpages{619--627}.
\bid{doi={10.2307/3315721}, issn={0319-5724}, mr={1671976}}
\bptok{imsref}%
\end{barticle}
\endbibitem

%b29 ###
\bibitem[\protect\citeauthoryear{Zhang, Kim and Woodroofe}{2001}]{Zhang2001}
\begin{barticle}[mr]
\bauthor{\bsnm{Zhang},~\bfnm{Rong}\binits{R.}},
  \bauthor{\bsnm{Kim},~\bfnm{Jeankyung}\binits{J.}} \AND
  \bauthor{\bsnm{Woodroofe},~\bfnm{Michael}\binits{M.}}
(\byear{2001}).
\btitle{Asymptotic analysis of isotonic estimation for grouped data}.
\bjournal{J. Statist. Plann. Inference}
\bvolume{98}
\bpages{107--117}.
\bid{doi={10.1016/S0378-3758(00)00301-3}, issn={0378-3758}, mr={1860229}}
\bptok{imsref}%
\end{barticle}
\endbibitem

\end{thebibliography}
\end{document}